\documentclass[reqno]{amsart}

\usepackage{tikz-cd}

\usepackage{fullpage,dsfont}
\usepackage{leftindex}

\usepackage{hyperref} 

\usepackage{parskip}
\makeatletter 
\def\thm@space@setup{%
 \thm@preskip=\parskip \thm@postskip=0pt
}
\def\th@remark{%
  \thm@headfont{\itshape}%
  \normalfont 
  \thm@preskip\parskip \thm@postskip=0pt
}
\makeatother

\usepackage[nobysame,alphabetic,initials,msc-links]{amsrefs}

\usepackage[final]{pdfpages}

\usepackage{multirow}

\usepackage{stmaryrd}

\usepackage{array}

\DefineSimpleKey{bib}{how}
\DefineSimpleKey{bib}{mrclass}
\DefineSimpleKey{bib}{mrnumber}
\DefineSimpleKey{bib}{fjournal}
\DefineSimpleKey{bib}{mrreviewer}

\renewcommand{\PrintDOI}[1]{%
  \href{http://dx.doi.org/#1}{{\tt DOI:#1}}%
}
\renewcommand{\eprint}[1]{#1}
\BibSpec{book}{%
    +{}  {\PrintPrimary}                {transition}
    +{.} { \PrintDate}                  {date}
    +{.} { \textit}                     {title}
    +{.} { }                            {part}
    +{:} { \textit}                     {subtitle}
    +{,} { \PrintEdition}               {edition}
    +{}  { \PrintEditorsB}              {editor}
    +{,} { \PrintTranslatorsC}          {translator}
    +{,} { \PrintContributions}         {contribution}
    +{,} { }                            {series}
    +{,} { \voltext}                    {volume}
    +{,} { }                            {publisher}
    +{,} { }                            {organization}
    +{,} { }                            {address}
    +{,} { }                            {status}
    +{,} { \PrintDOI}                   {doi}
    +{,} { \PrintISBNs}                 {isbn}
    +{}  { \parenthesize}               {language}
    +{}  { \PrintTranslation}           {translation}
    +{;} { \PrintReprint}               {reprint}
    +{.} { }                            {note}
    +{.} {}                             {transition}
    +{}  {\SentenceSpace \PrintReviews} {review}
}
\BibSpec{article}{%
    +{}  {\PrintAuthors}                {author}
    +{,} { \textit}                     {title}
    +{.} { }                            {part}
    +{:} { \textit}                     {subtitle}
    +{,} { \PrintContributions}         {contribution}
    +{.} { \PrintPartials}              {partial}
    +{,} { }                            {journal}
    +{}  { \textbf}                     {volume}
    +{}  { \PrintDatePV}                {date}
    +{,} { \issuetext}                  {number}
    +{,} { \eprintpages}                {pages}
    +{,} { }                            {status}
    +{,} { \PrintDOI}                   {doi}
    +{,} { \eprint}        {eprint}
    +{}  { \parenthesize}               {language}
    +{}  { \PrintTranslation}           {translation}
    +{;} { \PrintReprint}               {reprint}
    +{.} { }                            {note}
    +{.} {}                             {transition}
    +{}  {\SentenceSpace \PrintReviews} {review}
}
\BibSpec{collection.article}{%
    +{}  {\PrintAuthors}                {author}
    +{,} { \textit}                     {title}
    +{.} { }                            {part}
    +{:} { \textit}                     {subtitle}
    +{,} { \PrintContributions}         {contribution}
    +{,} { \PrintConference}            {conference}
    +{}  {\PrintBook}                   {book}
    +{,} { }                            {booktitle}
    +{,} { \PrintDateB}                 {date}
    +{,} { pp.~}                        {pages}
    +{,} { }                            {publisher}
    +{,} { }                            {organization}
    +{,} { }                            {address}
    +{,} { }                            {status}
    +{,} { \PrintDOI}                   {doi}
    +{,} { \eprint}        {eprint}
    +{}  { \parenthesize}               {language}
    +{}  { \PrintTranslation}           {translation}
    +{;} { \PrintReprint}               {reprint}
    +{.} { }                            {note}
    +{.} {}                             {transition}
    +{}  {\SentenceSpace \PrintReviews} {review}
}
\BibSpec{misc}{%
  +{}{\PrintAuthors}  {author}
  +{,}{ \textit}      {title}
  +{.}{ }             {how}
  +{}{ \parenthesize} {date}
  +{,} { available at \eprint}        {eprint}
  +{,}{ available at \url}{url}
  +{,}{ }             {note}
  +{.}{}              {transition}
}
\usepackage{amssymb, amsfonts, amsxtra, amsmath}
\usepackage{mathrsfs}
\usepackage{mathdots}
\usepackage{wasysym}
\usepackage[all]{xy}
\usepackage{bbm}
\usepackage{calc}
\usepackage{accents}

\usepackage[bbgreekl]{mathbbol}			

\usepackage{stackengine}

\numberwithin{equation}{section}

\DeclareSymbolFontAlphabet{\mathbb}{AMSb}	
\DeclareSymbolFontAlphabet{\mathbbl}{bbold}	

\newtheorem{Theorem}{Theorem}[section]
\newtheorem*{Theorem*}{Theorem}
\newtheorem{Def}[Theorem]{Definition}
\newtheorem*{Def*}{Def}
\newtheorem{Lem}[Theorem]{Lemma}
\newtheorem{Prop}[Theorem]{Proposition}

\newtheorem{Rem}[Theorem]{Remark}

\newtheorem{Exa}[Theorem]{Example}

\newcommand\bp{\begin{proof}}
\newcommand\ep{\end{proof}}

\mathchardef\mhyph="2D

\DeclareMathOperator{\id}{\mathrm{id}}

\DeclareMathOperator{\Rep}{\mathrm{Rep}}

\DeclareMathOperator{\Irr}{\mathrm{Irr}}

\newcommand{\msN}{\mathscr{N}}

\newcommand{\C}{\mathbb{C}}
\newcommand{\G}{\mathbb{G}}

\newcommand{\Corr}{\mathrm{Corr}}

\newcommand{\ovot}{\bar{\otimes}}

\begin{document}

\title{A categorical interpretation of Morita equivalence for dynamical von Neumann algebras}
\author{Joeri De Ro}
\address{Vrije Universiteit Brussel}
\email{joeri.ludo.de.ro@vub.be}

\begin{abstract} Let $\G$ be a locally compact quantum group and $(M, \alpha)$ a $\G$-$W^*$-algebra. The object of study of this paper is the $W^*$-category $\Rep^\G(M)$ of normal, unital $\G$-representations of $M$ on Hilbert spaces endowed with a unitary $\G$-representation. This category has a right action of the category $\Rep(\G)= \Rep^\G(\mathbb{C})$ for which it becomes a right $\Rep(\G)$-module $W^*$-category. Given another $\G$-$W^*$-algebra $(N, \beta)$, we denote the category of normal $*$-functors $\Rep^\G(N)\to \Rep^\G(M)$ compatible with the $\Rep(\G)$-module structure by $\operatorname{Fun}_{\Rep(\G)}(\Rep^\G(N), \Rep^\G(M))$ and we denote the category of $\G$-$M$-$N$-correspondences, studied in \cite{DCDR24}, by $\operatorname{Corr}^\G(M,N)$. We prove that there are canonical functors $P: \Corr^\G(M,N)\to \operatorname{Fun}_{\Rep(\G)}(\Rep^\G(N), \Rep^\G(M))$ and $Q: \operatorname{Fun}_{\Rep(\G)}(\Rep^\G(N), \Rep^\G(M))\to \operatorname{Corr}^\G(M,N)$ such that $Q \circ P\cong \id.$ We use these functors to show that the $\G$-dynamical von Neumann algebras $(M, \alpha)$ and $(N, \beta)$ are equivariantly Morita equivalent if and only if $\Rep^\G(N)$ and $\Rep^\G(M)$ are equivalent as $\Rep(\G)$-module-$W^*$-categories. Specializing to the case where $\G$ is a compact quantum group, we prove that moreover $P\circ Q \cong \id$, so that the categories $\Corr^\G(M,N)$ and $\operatorname{Fun}_{\Rep(\G)}(\Rep^\G(N), \Rep^\G(M))$ are equivalent. This is an equivariant version of the Eilenberg-Watts theorem for actions of compact quantum groups on von Neumann algebras.
\end{abstract}

\maketitle

\section{Introduction}

Given von Neumann algebras $M$ and $N$, we write $\operatorname{Corr}(M,N)$ for the category of $M$-$N$-correspondences, that is Hilbert spaces $\mathcal{H}$ with a unital normal $*$-representation $\pi: M\to B(\mathcal{H})$ and a unital normal anti-$*$-representation $\rho: N \to B(\mathcal{H})$ such that $\pi(m) \rho(n) = \rho(n)\pi(m)$ for all $m\in M$ and all $n\in N$. Let us write $\Rep(M) = \operatorname{Corr}(M,\mathbb{C})$ for the $W^*$-category of unital, normal $*$-representations of $M$ on Hilbert spaces. This category has been systematically studied in \cite{Rie74}. In particular, it is shown that normal $*$-functors $F: \Rep(N)\to \Rep(M)$ can be described by the notion of \emph{self-dual $N$-rigged $M$-modules}. Since the category of self-dual $N$-rigged $M$-modules is equivalent with the category $\operatorname{Corr}(M,N)$ \cite{BDH88}*{Théorème 2.2}, we can equivalently describe $*$-functors $F: \Rep(N)\to \Rep(M)$ by $M$-$N$-correspondences. This equivalent viewpoint is used in \cites{Bro03, Sau83} and will also be convenient for us. We now give a brief discussion about this.

If $\mathcal{G}\in \operatorname{Corr}(M,N)$, the Connes fusion tensor product induces a normal $*$-functor
$$F_{\mathcal{G}}: \Rep(N)\to \Rep(M): \mathcal{H}\mapsto \mathcal{G}\boxtimes_N \mathcal{H}.$$
Conversely, if $F: \Rep(N)\to \Rep(M)$ is a normal $*$-functor, then the space $\mathcal{G}:=F(L^2(N))\in \Rep(M)$ carries the unital, normal anti-$*$-representation
$$\rho: N \to B(\mathcal{G}): n \mapsto F(\rho_N(n))$$
so that $\mathcal{G}\in \Corr(M,N)$. Moreover, it can be shown that $F$ is unitarily naturally isomorphic to $F_{\mathcal{G}}$, establishing a von Neumann algebra version of the celebrated Eilenberg-Watts theorem \cites{Ei60, Wa60}. Under the above constructions, a Morita equivalence Hilbert space $\mathcal{G}$ corresponds to an equivalence $\Rep(M)\sim \Rep(N)$ of $W^*$-categories. Consequently, the von Neumann algebras $M$ and $N$ are Morita equivalent if and only if $\Rep(M)$ and $\Rep(N)$ are equivalent as $W^*$-categories \cite{Bro03}*{Theorem 3.5.5}. Denoting the category of normal $*$-functors $\Rep(N)\to \Rep(M)$  by $\operatorname{Fun}(\Rep(N), \Rep(M))$, the above constructions induce functors
$$P: \Corr(M,N)\to \operatorname{Fun}(\Rep(N), \Rep(M)), \quad Q: \operatorname{Fun}(\Rep(N), \Rep(M))\to \Corr(M,N)$$
that are quasi-inverse to each other, so that the categories $\Corr(M,N)$ and $\operatorname{Fun}(\Rep(N), \Rep(M))$ are equivalent to each other.

The main goal of this paper is to examine how these results generalise to the \emph{equivariant} setting. More concretely, given a locally compact quantum group $\G$, a $\G$-$W^*$-algebra $(M, \alpha)$ consists of a von Neumann algebra $M$ together with a unital, isometric, normal $*$-homomorphism $\alpha: M\to M \ovot L^\infty(\G)$ such that $(\alpha\otimes \id)\alpha = (\id \otimes \Delta)\alpha$. We write $\alpha: M\curvearrowleft \G$ and we say that $\alpha$ is an action of $\G$ on $M$. Given another $\G$-$W^*$-algebra $(N, \beta)$, we then consider the category of $\G$-$M$-$N$-correspondences $\operatorname{Corr}^\G(M,N)$. Objects in this category consist of correspondences $\mathcal{H}= (\mathcal{H}, \pi, \rho)\in \operatorname{Corr}(M,N)$ together with a unitary $\G$-representation $U\in B(\mathcal{H})\ovot L^\infty(\G)$ such that 
$$(\pi\otimes \id)\alpha(m) = U(\pi(m)\otimes 1)U^*, \quad (\rho\otimes R)\beta(n) = U^*(\rho(n)\otimes 1)U, \quad m\in M, \quad n \in N,$$
and where $R: L^\infty(\G)\to L^\infty(\G)$ is the unitary antipode. We then write $\Rep^\G(M)= \operatorname{Corr}^\G(M, \C)$. Again, if $\mathcal{G}\in \Corr^\G(M,N)$, there is an induced normal $*$-functor
$$F_{\mathcal{G}}: \Rep^\G(N)\to \Rep^\G(M): \mathcal{H}\mapsto \mathcal{G}\boxtimes_N \mathcal{H}.$$
It is natural to ask if every normal $*$-functor $F: \Rep^\G(N)\to \Rep^\G(M)$ then induces a correspondence $\mathcal{G}\in \Corr^\G(M,N)$. Unfortunately, this is no longer true, but it is still possible to find a suitable generalization by taking into consideration some extra structure that is not present (or rather: trivial) in the non-equivariant case. More concretely, consider the $W^*$-category $\Rep(\G)=\Rep^\G(\mathbb{C})= \operatorname{Corr}^\G(\C, \C)$ of unitary $\G$-representations, and note that $\Rep^\G(M)$ becomes a right $\Rep(\G)$-module $W^*$-category in a natural way (by tensoring). The functor $F_{\mathcal{G}}$ then satisfies
$$F_{\mathcal{G}}(\mathcal{H}\otimes \mathcal{K})\cong F_{\mathcal{G}}(\mathcal{H})\otimes \mathcal{K}$$
by unitaries natural in $\mathcal{H}\in \Rep^\G(N)$ and natural in $\mathcal{K}\in \Rep(\G)$. In other words, $F_{\mathcal{G}}$ is a functor compatible with the $\Rep(\G)$-module structure. Let us denote the category of normal $*$-functors $\Rep^\G(N)\to \Rep^\G(M)$ that are compatible with the $\Rep(\G)$-module structure by $\operatorname{Fun}_{\Rep(\G)}(\Rep^\G(N), \Rep^\G(M)).$ Such functors induce $\G$-$M$-$N$-correspondences, as the following result shows:

\textbf{Theorem \ref{mr} + Proposition \ref{functors}} \textit{ 
Let $\G$ be a locally compact quantum group, let $M,N$ be $\G$-$W^*$-algebras and let $F\in \operatorname{Fun}_{\Rep(\G)}(\Rep^\G(N), \Rep^\G(M))$. Consider $(\mathcal{G}, \pi_{\mathcal{G}}, U_{\mathcal{G}}):=  F(L^2(N))\in \Rep^\G(M)$. The Hilbert space $\mathcal{G}$ admits a canonical anti-$*$-representation
    $\rho_{\mathcal{G}}: N \to B(\mathcal{G})$
    such that $(\mathcal{G}, \pi_{\mathcal{G}}, \rho_{\mathcal{G}}, U_{\mathcal{G}})\in \Corr^\G(M,N).$ In this way, we obtain canonical functors}
\begin{align*}
    &P: \Corr^\G(M,N)\to \operatorname{Fun}_{\Rep(\G)}(\Rep^\G(N), \Rep^\G(M)),\\
    &Q: \operatorname{Fun}_{\Rep(\G)}(\Rep^\G(N), \Rep^\G(M))\to \Corr^\G(M,N).
\end{align*}
\textit{Moreover, $Q\circ P \cong \id.$}

Note that it is no longer clear how the right action $\mathcal{G}\curvearrowleft N$ has to be constructed in the equivariant setting. Our strategy is to first construct a right action $\mathcal{G}\otimes L^2(\G)\curvearrowleft 
 N\rtimes_\beta \G$ which is then shown to be induced by the desired right action $\mathcal{G}\curvearrowleft N$.

 We can then prove the following categorical characterization of equivariant Morita equivalence:

 \textbf{Theorem \ref{equivariant_morita}} \textit{Let $\G$ be a locally compact quantum group. The $\G$-$W^*$-algebras $(M, \alpha)$ and $(N, \beta)$ are equivariantly $W^*$-Morita equivalent if and only if $\Rep^\G(M)$ and $\Rep^\G(N)$ are equivalent as $\Rep(\G)$-module $W^*$-categories.}

In algebraic contexts, it was already known that equivariant Morita equivalence of certain objects could be detected by considering associated module categories, see e.g.\ \cite{AM07}.

A natural question is if it is also true that  $P\circ Q\cong \id$, and thus that every normal $*$-functor $F: \Rep^\G(N)\to \Rep^\G(M)$ compatible with the $\Rep(\G)$-module stucture arises by tensoring with a $\G$-$M$-$N$-correspondence. This is unlikely to be true in general, although we do not know a concrete counterexample. However, for compact quantum groups this becomes true:

\textbf{Theorem \ref{main}}  \textit{Let $\G$ be a compact quantum group. The functors}
    \begin{align*}
        &P: \operatorname{Corr}^\G(M,N)\to \operatorname{Fun}_{\Rep(\G)}(\Rep^\G(N), \Rep^\G(M)): \mathcal{G}\mapsto F_{\mathcal{G}}, \\
    &Q: \operatorname{Fun}_{\Rep(\G)}(\Rep^\G(N), \Rep^\G(M))\to  \operatorname{Corr}^\G(M,N): F \mapsto F(L^2(N))
    \end{align*}
    \textit{are quasi-inverse to each other. In particular, every normal $*$-functor $F: \Rep^\G(N)\to \Rep^\G(M)$ that is compatible with the $\Rep(\G)$-module structure is naturally isomorphic to $F_{\mathcal{G}}$ where $\mathcal{G}= F(L^2(N))\in \Corr^\G(M,N).$}

    This result can be seen as an equivariant Eilenberg-Watts theorem in the setting of actions of compact quantum groups on von Neumann algebras.

    \section{Preliminaries}
    We now recall some useful definitions and results.

\subsection{Locally compact quantum groups} Given a von Neumann algebra $M$, we consider the standard Hilbert space $L^2(M)$ with modular conjugation $J_M$. It carries a unital, normal $*$-representation $\pi_M: M\to B(L^2(M))$ and a unital, normal, anti-$*$-representation $\rho_M: M\to B(L^2(M)): x \mapsto J_M \pi_M(x)^* J_M$. We have $\rho_M(M) = \pi_M(M)'$. We recall some theory of locally compact quantum groups \cites{KV00,KV03,VV03}.

A \emph{von Neumann bialgebra} is a pair $(M, \Delta)$ where $M$ is a von Neumann algebra and $\Delta: M \to M \ovot M$ a unital, normal, isometric $*$-homomorphism such that $(\Delta \otimes \id)\circ \Delta = (\id \otimes \Delta)\circ \Delta$. Here, $M\ovot M$ denotes the von Neumann algebra tensor product, containing the algebraic tensor product $M \odot M$ as a $\sigma$-weakly dense $*$-subalgebra. 

A \emph{locally compact quantum group} $\G$ is a von Neumann bialgebra $(L^\infty(\G), \Delta)$ for which there exist normal, semifinite, faithful weights $\Phi, \Psi: L^\infty(\G)_+\to [0, \infty]$ such that $(\id\otimes \Phi)\Delta(x)= \Phi(x)1$ for all $x\in \mathscr{M}_\Phi^+$ and $(\Psi\otimes \id)\Delta(x)= \Psi(x)1$ for all $x\in \mathscr{M}_\Psi^+$. The von Neumann algebra predual $L^1(\G):= L^\infty(\G)_*$ is a completely contractive Banach algebra for the multiplication 
$$\mu\star \nu:= (\mu\otimes \nu)\circ \Delta, \quad \mu, \nu \in L^1(\G).$$

Without any loss of generality, we may (and will) assume that $L^\infty(\G)$ is faithfully represented in standard form on the Hilbert space $L^2(\G)$. With respect to the weights $\Phi, \Psi$, we then have GNS-maps
$$\Lambda_\Phi: \mathscr{N}_\Phi\to L^2(\G), \quad \Lambda_\Psi: \mathscr{N}_\Psi\to L^2(\G).$$

Fundamental to the theory of locally compact quantum groups are the unitaries 
\[
V,W \in B(L^2(\G)\otimes L^2(\G))
\]
called respectively \emph{right} and \emph{left} regular unitary representation. They are uniquely characterized by the identities 
\[
(\id\otimes \omega)(V) \Lambda_{\Psi}(x) = \Lambda_{\Psi}((\id\otimes \omega)\Delta(x)),
\qquad \omega \in L^1(\G),x\in \msN_{\Psi},
\]
\[
 (\omega \otimes \id)(W^*)\Lambda_{\Phi}(x) = \Lambda_{\Phi}((\omega\otimes \id)\Delta(x)),\qquad \omega \in L^1(\G),x\in \msN_{\Phi}.
\]
They are \emph{multiplicative unitaries} \cite{BS93} meaning that
\[
V_{12}V_{13}V_{23} = V_{23}V_{12},\qquad W_{12}W_{13}W_{23}= W_{23}W_{12}
\]
and they implement the coproduct of $L^\infty(\G)$ in the sense that
$$\label{EqComultImpl}
W^*(1\otimes x)W = \Delta(x) = V(x\otimes 1)V^*,\qquad x\in L^\infty(\G).$$
Moreover, we have
$$C_0(\G):=[(\omega\otimes \id)(V) \mid \omega \in B(L^2(\G))_*] = [(\id\otimes \omega)(W) \mid \omega \in B(L^2(\G))_*]$$
which is a $\sigma$-weakly dense C$^*$-subalgebra of $L^\infty(\G)$. Then $\Delta(C_0(\G))\subseteq M(C_0(\G)\otimes C_0(\G))$. 

We can also define the von Neumann algebra
$L^\infty(\hat{\G}) = [(\omega\otimes \id)(W) \mid \omega \in L^1(\G)]^{\sigma\textrm{-weak}}$
with coproduct
$\hat{\Delta}(x) = \Sigma W(x\otimes 1)W^*\Sigma$ where $x\in L^\infty(\hat{\G})$. 
Then the pair $(L^\infty(\hat{\G}), \hat{\Delta})$ defines the locally compact quantum group $\hat{\G}$, which is called the dual of $\G$. Similarly, we can also define the von Neumann algebra $L^\infty(\check{\G}) = [(\id\otimes \omega)(V) \mid \omega \in L^1(\G)]^{\sigma\textrm{-weak}}$
with coproduct 
$\check{\Delta}(x) = V^*(1\otimes x)V$ for $x\in L^\infty(\check{\G})$.
Then the pair $(L^\infty(\check{\G}), \check{\Delta})$ defines a locally compact quantum group $\check{\G}$. We have $L^\infty(\hat{\G})= L^\infty(\check{\G})'$. If $G$ is a locally compact group, we have
$L^\infty(\hat{G}) = \mathscr{L}(G)$ and $L^\infty(\check{G}) = \mathscr{R}(G)$, the left resp. right group von Neumann algebra associated with $G$. Therefore, $\hat{\G}$ and $\check{\G}$ should be thought of as a left and a right version of the dual locally compact quantum group.

The left multiplicative unitary of $\check{\G}$ will be denoted by $\check{W}$ and the right multiplicative unitary of $\check{\G}$ will be denoted by $\check{V}$. We have $\check{W}= V$ and 
$$V \in L^\infty(\check{\G})\ovot L^\infty(\G), \quad W \in L^\infty(\G)\ovot L^\infty(\hat{\G}), \quad \check{V}\in L^\infty(\G)'\ovot L^\infty(\check{\G}).$$
In fact, we have e.g.\ $V\in M(C_0(\check{\G})\otimes C_0(\G))$ and similarly for $W$ and  $\check{V}$. For future use, we will also need the following notations:
\begin{align*}
    &\Delta_l: B(L^2(\G)) \to  L^\infty(\G)\ovot B(L^2(\G)): x \mapsto W^*(1\otimes x)W,\\
    &\Delta_r: B(L^2(\G))\to B(L^2(\G))\ovot L^\infty(\G): x \mapsto V(x\otimes 1)V^*,\\
    &\check{\Delta}_r: B(L^2(\G))\to B(L^2(\G)) \ovot L^\infty(\check{\G}): x \mapsto \check{V}(x\otimes 1)\check{V}^*.
\end{align*}
A locally compact quantum group $\G$ is called \emph{compact} if $C_0(\G)$ is unital, and we then write $C_0(\G)= C(\G)$. In that case, the left Haar weight $\Phi$ is a normal state and $\Phi = \Psi$ (after appropriate normalisations)

A locally compact quantum group $\G$ is called \emph{discrete} if $\check{\G}$ is compact and we then write $L^\infty(\G)= \ell^\infty(\G)$ and $C_0(\G)= c_0(\G)$. If $\G$ is a discrete quantum group, there is a unique normal state $\epsilon \in \ell^1(\G):= L^1(\G)$, called \emph{counit}, such that 
$(\epsilon \otimes \id)\Delta = \id = (\id \otimes \epsilon)\Delta$.

A \emph{unitary representation} of a locally compact quantum group $\G$ on a Hilbert space consists of a pair $(\mathcal{H}, U)$ where $\mathcal{H}$ is a Hilbert space and $U\in B(\mathcal{H})\ovot L^\infty(\G)$ a unitary satisfying $(\id \otimes \Delta)(U) = U_{12}U_{13}$. We write $(\mathcal{H}, U)\in \Rep(\G)$. For example, we have $V \in \Rep(\G)$ and  $W_{21}\in \Rep(\G)$, where $W_{21}= \sigma W$ with $\sigma: B(L^2(\G))\ovot B(L^2(\G))\cong B(L^2(\G))\ovot B(L^2(\G))$ the flip. Given a Hilbert space $\mathcal{H}$, we write $\mathbb{I}= 1 \otimes 1 \in B(\mathcal{H})\ovot L^\infty(\G)$ for the trivial representation. Given $(\mathcal{H}, U_{\mathcal{H}}), (\mathcal{K}, U_{\mathcal{K}})\in \Rep(\G)$, we also have that $(\mathcal{H}\otimes \mathcal{K}, U_{\mathcal{H}, 13}U_{\mathcal{K}, 23})\in \Rep(\G)$. Given $U \in B(\mathcal{H})\ovot L^\infty(\G)$ and $\omega \in L^1(\G)$, we write $U(\omega):= (\id \otimes \omega)(U)$.

Consider the modular conjugation $J:= J_{L^\infty(\G)}$ associated to the standard representation $L^\infty(\G)\subseteq B(L^2(\G))$ and the modular conjugation $\check{J}:= J_{L^\infty(\check{\G})}$ associated to the standard representation $L^\infty(\check{\G})\subseteq B(L^2(\G))$. We have 
$\check{J}L^\infty(\G)\check{J}=L^\infty(\G)$
so we obtain the anti-$*$-morphism
$$R: L^\infty(\G)\to L^\infty(\G): x \mapsto \check{J}x^*\check{J}.$$
We call $R$ the \emph{unitary antipode} of $\G$. There is a canonical unimodular complex number $c\in \mathbb{C}$ such that 
$$c \check{J}J=\overline{c}J \check{J}.$$
We write $u_\G:= c \check{J} J$ for the associated self-adjoint unitary. The following identities will be useful:
\begin{align*}
    (\check{J}\otimes J)W(\check{J}\otimes J)=W^*, \quad (J\otimes \check{J})V(J\otimes \check{J})=V^*, \quad (u_\G\otimes 1)V(u_\G\otimes 1)=W_{21}.
\end{align*}

\subsection{Dynamical von Neumann algebras and crossed products}
Let $\G$ be a locally compact quantum group.
A (right) \emph{$\G$-$W^*$-algebra} is a pair $(M, \alpha)$ such that $M$ is a von Neumann algebra and $\alpha: M \to M\ovot L^\infty(\G)$ is an injective, unital, normal $*$-homomorphism satisfying the coaction property $(\alpha\otimes \id)\circ \alpha = (\id \otimes \Delta)\circ \alpha$. We sometimes denote this with $\alpha: M\curvearrowleft\G$. We write
$M^{\alpha} = \{x\in M: \alpha(x)= x \otimes 1\}$
for the von Neumann subalgebra of fixed points of $(M, \alpha)$ and the trivial $\G$-action on a von Neumann algebra $M$ will always be denoted by $\tau$, i.e. $\tau(x)= x \otimes 1$ for $x\in M$. Given a $\G$-$W^*$-algebra $(M, \alpha)$,  we define the \emph{crossed product von Neumann algebra}
$M\rtimes_\alpha \G= [\alpha(M)(1\otimes L^\infty(\check{\G}))]''.$
We have the alternative description
$$M\rtimes_\alpha \G = \{z \in M \ovot B(L^2(\G)): (\alpha\otimes \id)(z) =(\id \otimes \Delta_l)(z)\}.$$
It is easy to verify that 
$$(\id \otimes \check{\Delta}_r)(M\rtimes_\alpha \G)\subseteq (M\rtimes_\alpha \G)\ovot L^\infty(\check{\G})$$
so that $(M\rtimes_\alpha \G, \id \otimes \check{\Delta}_r)$ becomes a $\check{\G}$-$W^*$-algebra. This $\check{\G}$-action on the crossed product will often be implicitly understood. We have that $(M\rtimes_\alpha \G)^{\id \otimes \check{\Delta}_r}= \alpha(M)$. Given any von Neumann algebra $M$, we have $M\rtimes_\tau \G = M\ovot L^\infty(\check{\G})$.

If $(M, \alpha)$ is a $\G$-$W^*$-algebra where $M$ is standardly represented on a Hilbert space $\mathcal{H}$, there exists a canonical unitary $U_\alpha \in B(\mathcal{H})\ovot L^\infty(\G)$ such that $\alpha(m)= U_\alpha(m\otimes 1)U_\alpha^*$ for all $m\in M$ \cite{Va01}. We call $U_\alpha$ the \emph{unitary implementation} of $\alpha$, and it satisfies $(\id \otimes \Delta)(U_\alpha)= U_{\alpha,12}U_{\alpha,13}.$ More precisely, we can identify $L^2(M\rtimes_\alpha \G)= L^2(M)\otimes L^2(\G)$ in such a way that the standard representation $\pi_{M\rtimes_\alpha \G}$ corresponds to $\pi_M\otimes \id: M\rtimes_\alpha \G \to B(L^2(M)\otimes L^2(\G))$ and such that the anti-$*$-representation $\rho_{M\rtimes_\alpha \G}: M\rtimes_\alpha \G \to B(L^2(M)\otimes L^2(\G))$ is given by
$$\rho_{M\rtimes_\alpha \G}(\alpha(m)) = \rho_M(m)\otimes 1, \quad \rho_{M\rtimes_\alpha \G}(1\otimes \check{y})=U_\alpha(1\otimes \rho_{L^\infty(\check{\G})}(\check{y}))U_\alpha^*, \quad m\in M, \quad \check{y}\in L^\infty(\check{\G}).$$

\subsection{Equivariant correspondences} Let $\G$ be a locally compact quantum group.

Let $(M, \alpha)$ and $(N, \beta)$ be two $\G$-$W^*$-algebras. 
A $\G$-$M$-$N$-correspondence consists of a Hilbert space $\mathcal{H}$ together with the following data:
\begin{itemize}
    \item A unitary $\G$-representation $U\in B(\mathcal{H})\ovot L^\infty(\G)$.
    \item A unital, normal $*$-homomorphism $\pi: M\to B(\mathcal{H})$ such that $(\pi\otimes \id)\alpha(m) = U(\pi(m)\otimes 1)U^*$ for all $m\in M$.
    \item A unital, normal anti-$*$-homomorphism $\rho: N \to B(\mathcal{H})$ such that $(\rho\otimes R)\beta(n) = U^*(\rho(n)\otimes 1)U$ for all $n\in N$.
\end{itemize}

Given two $\G$-$M$-$N$-correspondences $\mathcal{H}$ and $\mathcal{K}$, a bounded linear map $x: \mathcal{H}\to \mathcal{K}$ is called intertwiner of correspondences if 
$$x\pi_{\mathcal{H}}(m) = \pi_{\mathcal{K}}(m)x, \quad x\rho_{\mathcal{H}}(n) = \rho_{\mathcal{K}}(n)x, \quad x U_{\mathcal{H}}(\omega) = U_{\mathcal{K}}(\omega)x, \quad m\in M, \quad n\in N, \quad \omega \in L^1(\G).$$
The set of such intertwiners is denoted by ${}_M\mathscr{L}_N^\G(\mathcal{H}, \mathcal{K})$. If $\G$ is the trivial group or if one of the von Neumann algebras $M$ or $N$ is equal to $\C$, we omit it from the notation. 

In this way, we obtain the $W^*$-category $\Corr^\G(M,N)$ of $\G$-$M$-$N$-correspondences. 

We call $(M, \alpha)$ and $(N, \beta)$ \emph{$\G$-$W^*$-Morita equivalent}, and we write $(M,\alpha)\sim_\G (N, \beta)$, if there exists $(\mathcal{H}, \pi, \rho, U)\in \operatorname{Corr}^\G(M,N)$ such that $\pi$ and $\rho$ are faithful (= injective) and such that $\pi(M)' = \rho(N).$ If $M,N$ are Morita equivalent von Neumann algebras, then we simply write $M\sim N$. This is the situation where $\G$ is the trivial group. It is not so hard to see that $\sim_\G$ is an equivalence relation. The notion of equivariant Morita equivalence for dynamical von Neumann algebras will be studied in depth in future work.

There are several examples and constructions we can perform with equivariant correspondences. We give a brief overview of those relevant to this paper:

\textbf{Direct sums.} Given a collection $\{\mathcal{H}_i\}_{i\in I}$ of $\G$-$M$-$N$-correspondences, the Hilbert space direct sum $\mathcal{H}:=\bigoplus_{i\in I}\mathcal{H}_i$ is endowed with the unique $\G$-$M$-$N$-correspondence structure such that the inclusion maps $\mathcal{H}_i\hookrightarrow \mathcal{H}$ all become intertwiners of $\G$-$M$-$N$-correspondences.

\textbf{Crossed products.} Given $\mathcal{H}=(\mathcal{H}, \pi, \rho, U)\in \Corr^\G(M,N)$, it was shown in \cite[Proposition 5.21]{DCDR24} that the following data defines a $\check{\G}$-correspondence $\mathcal{H}^\rtimes \in \Corr^{\check{\G}}(M\rtimes_\alpha \G, N\rtimes_\beta \G):$
\begin{itemize}
    \item The Hilbert space $\mathcal{H}\otimes L^2(\G)$.
    \item The unitary $\check{\G}$-representation $\check{V}_{23}\in B(\mathcal{H}\otimes L^2(\G))\ovot L^\infty(\check{\G})$.
    \item The normal $*$-representation $\pi^\rtimes: M\rtimes_\alpha \G \to B(\mathcal{H}\otimes L^2(\G)): z \mapsto (\pi\otimes \id)(z)$.
    \item The normal $*$-anti-representation $\rho^\rtimes: N\rtimes_\beta \G \to B(\mathcal{H}\otimes L^2(\G)): z \mapsto U(\rho\otimes \check{J}(-)^*\check{J})(z)U^*$.
\end{itemize}

For future use, we will need to know how the crossed product construction behaves with respect to equivariant Morita equivalence:

\begin{Prop}\label{crossed product}
    Let $\mathcal{H}= (\mathcal{H}, \pi, \rho, U)\in \operatorname{Corr}^\G(M,N)$. The following are equivalent:
    \begin{enumerate}
        \item $\mathcal{H}$ is a $\G$-$W^*$-Morita equivalence $(M, \alpha)\sim_\G (N, \beta)$
        \item $\mathcal{H}^\rtimes$ is a $\check{\G}$-$W^*$-Morita equivalence $(M\rtimes_\alpha \G, \id \otimes \check{\Delta}_r)\sim_{\check{\G}}(N\rtimes_\beta \G, \id \otimes \check{\Delta}_r).$
    \end{enumerate}
    Therefore, $(M, \alpha)\sim_\G (N, \beta)\iff (M\rtimes_\alpha \G, \id \otimes \check{\Delta}_r)\sim_{\check{\G}}(N\rtimes_\beta \G, \id \otimes \check{\Delta}_r).$
\end{Prop}
\begin{proof} $(1)\implies (2)$ Assume first that $\mathcal{H}$ is a $\G$-$W^*$-Morita equivalence. Assume that $z\in \rho^{\rtimes}(N\rtimes_\beta \G)'$. Then we have
\begin{align}
    \label{1}&z (\rho(n)\otimes 1) = (\rho(n)\otimes 1)z, \quad n \in N,\\
    \label{2}&zU(1\otimes \hat{y})U^* = U(1\otimes \hat{y})U^*z, \quad  \hat{y}\in L^\infty(\hat{\G}).
\end{align}
In particular, \eqref{1} implies that 
$$z \in (\rho(N)\ovot \mathbb{C}1)' = \pi(M)\ovot B(L^2(\G))$$
and thus we can write $z= (\pi\otimes \id)(y)$ for some unique $y\in M\ovot B(L^2(\G))$. On the other hand, \eqref{2} implies that for all $\omega \in L^1(\G)$, we have
$$(\id \otimes \omega \otimes \id)(z_{13} U_{13}W_{23}U_{13}^*)= z U(1\otimes (\omega \otimes \id)(W))U^* = U(1\otimes (\omega \otimes \id)(W))U^* z = (\id \otimes \omega \otimes \id)(U_{13}W_{23}U_{13}^* z_{13})$$
whence 
$$z_{13}U_{13}W_{23}U_{13}^* = U_{13}W_{23}U_{13}^* z_{13}.$$
From $(\id \otimes \Delta)(U) = U_{12}U_{13}$, we find that
$W_{23}^* U_{13}W_{23}U_{13}^* = U_{12}$
so that 
\begin{align*}
    (\pi\otimes \id \otimes \id)(\alpha\otimes \id)(y)&= U_{12}(\pi\otimes \id)(y)_{13}U_{12}^*=W_{23}^* U_{13}W_{23}U_{13}^* z_{13} U_{13}W_{23}^* U_{13}^*W_{23}\\
    &= W_{23}^* z_{13}W_{23} = (\id \otimes \Delta_l)(z)= (\pi\otimes \id \otimes \id)(\id \otimes \Delta_l)(y)
\end{align*}
so that injectivity of $\pi$ implies that 
$$(\alpha\otimes \id)(y) = (\id \otimes \Delta_l)(y)$$
whence 
$y \in M\rtimes_\alpha \G$
and thus $z= \pi^\rtimes(y).$ Therefore, $\pi^\rtimes(M\rtimes \G)= \rho^\rtimes(N\rtimes_\beta \G)'.$

$(2)\implies (1)$ Assume that $\mathcal{H}^\rtimes$ is a $\check{\G}$-$W^*$-Morita equivalence. If $x\in \pi(M)'$, then $x\otimes 1 \in \pi^\rtimes(M\rtimes_\alpha\G)' = \rho^\rtimes(N\rtimes_\beta \G)$, so there is a unique $z\in N\rtimes_\beta \G$ such that $x\otimes 1 = \rho^\rtimes(z)$. Then
\begin{align*}
    (\rho^\rtimes \otimes \check{R})((\id \otimes \check{\Delta}_r)(z))&= \check{V}_{23}^*(\rho^{\rtimes}(z)\otimes 1)\check{V}_{23} = x\otimes 1 \otimes 1 = (\rho^{\rtimes}\otimes \check{R})(z\otimes 1)
\end{align*}
from which we conclude that $z\in \beta(N)$. Therefore, $z= \beta(n)$ and thus $x= \rho(n)$ for some $n\in N$. Therefore, $\pi(M)' = \rho(N)$. 

It is now clear that
$$(M, \alpha)\sim_\G (N, \beta)\implies (M\rtimes_\alpha \G, \id \otimes \check{\Delta}_r)\sim_{\check{\G}} (N\rtimes_\beta \G, \id \otimes \check{\Delta}_r).$$

Conversely, if $(M\rtimes_\alpha \G, \id \otimes \check{\Delta}_r)\sim_{\check{\G}} (N\rtimes_\beta \G, \id \otimes \check{\Delta}_r)$, then by the implication we just proved and Takesaki-Takai duality \cite{DCDR24}*{Theorem 5.25}, we find
$$M \sim_{\check{\check{\G}}}M\rtimes \G \times \check{\G}\sim_{\check{\check{\G}}} N\rtimes \G \rtimes \check{\G} \sim_{\check{\check{\G}}} N.$$
From the Pontryagin biduality $\G \cong \check{\check{\G}}$, we then conclude that $M \sim_\G N.$
\end{proof}

\textbf{Connes fusion product.} Given $\mathcal{H}\in \Corr^\G(M,N)$ and $\mathcal{G}\in \Corr^\G(N,P)$, there is a way to build a composite correspondence $\mathcal{H}\boxtimes_N \mathcal{G}\in \Corr^\G(M,P)$. More precisely, consider the algebraic balanced tensor product
$$\mathscr{L}_N(L^2(N), \mathcal{H})\odot_N \mathcal{G}$$
which we endow with the semi-inner-product uniquely determined by
$$\langle x\otimes_N \xi, y \otimes_N \eta\rangle:= \langle \xi, \pi_{\mathcal{G}}(\pi_N^{-1}(x^* y))\eta\rangle.$$
We then define the Hilbert space $\mathcal{H}\boxtimes_N \mathcal{G}$ to be the separation-completion of $\mathscr{L}_N(L^2(N), \mathcal{H})\odot_N \mathcal{G}$. It carries a left $M$-action via
$$\pi_{\boxtimes}(m)(x\otimes_N \xi)= \pi_{\mathcal{H}}(m)x\otimes_N \xi$$
and a right $P$-action via
$$\rho_{\boxtimes}(p)(x\otimes_N \xi) = x\otimes_N \rho_{\mathcal{G}}(p)\xi.$$
Moreover, this Hilbert space carries a canonical $\G$-representation $U_{\boxtimes}\in B(\mathcal{H}\boxtimes_N \mathcal{G})\ovot L^\infty(\G)$ such that $$(\mathcal{H}\boxtimes_N \mathcal{G}, \pi_\boxtimes, \rho_\boxtimes, U_\boxtimes)\in \Corr^\G(M,P).$$For more details, the reader is refered to \cite[Proposition 5.6]{DCDR24}.

\subsection{Categories of equivariant representations.} We introduce some more notation:
\begin{itemize}
    \item We write $\Rep(\G)= \Corr^\G(\mathbb{C}, \mathbb{C})$ for the category of unitary $\G$-representations.
    \item We write $\Rep^\G(M)= \Corr^\G(M, \C)$ for the category of triples $(\mathcal{H}, \pi, U)$ where $\mathcal{H}$ is a Hilbert space, $U\in B(\mathcal{H})\ovot L^\infty(\G)$ is a unitary $\G$-representation and $\pi: M\to B(\mathcal{H})$ is a unital, normal $*$-homomorphism such that $(\pi\otimes \id)\alpha(m) = U(\pi(m)\otimes 1)U^*$ for all $m\in M$. We also refer to these as $\G$-representations of $(M, \alpha)$.
\end{itemize}

Often, we will simply write $\mathcal{H}\in \Rep^\G(M)$ and the normal $*$-representation of $M$ on $\mathcal{H}$ and the unitary representation of $\G$ on $\mathcal{H}$ will be implicitly understood. In a similar fashion, unless explicitly mentioned otherwise, when we write $L^2(\G)\in \Rep(\G)$, we will always mean $(L^2(\G),V)$ where $V\in L^\infty(\check{\G})\ovot L^\infty(\G)$ is the right regular representation. 

Let us mention some important objects in $\Rep^\G(M)$:
\begin{enumerate}
    \item $(L^2(M), \pi_M, U_\alpha) \in \Rep^\G(M)$, the trivial $\G$-representation of $(M,\alpha)$.
    \item Given $\mathcal{G}\in \Rep^\G(M)$, we consider 
$$S^\G(\mathcal{G}):= (\mathcal{G}\otimes L^2(\G), \pi_S: M\to B(\mathcal{G}\otimes L^2(\G)): m \mapsto (\pi_\mathcal{G}\otimes \id)\alpha(m), V_{23}) \in \Rep^\G(M).$$
We write $S^\G_M:= S^\G(L^2(M))$ for the \emph{semi-coarse} $\G$-representation of $(M, \alpha)$.
\end{enumerate}

We will need the following easy lemma:

\begin{Lem}\label{easy}
If $\mathcal{G}\in \Rep^\G(M)$, we have    ${}_M\mathscr{L}^\G(S^\G(\mathcal{G}))= (\pi_\mathcal{G}\otimes \id)(M\rtimes_\alpha \G)'.$
\end{Lem}
\begin{proof} Note that $z\in (\pi_{\mathcal{G}}\otimes \id)(M\rtimes_\alpha \G)'$ if and only if 
\begin{align}\label{hah}
    &z(\pi_{\mathcal{G}}\otimes \id)(\alpha(m)) = (\pi_{\mathcal{G}}\otimes \id)(\alpha(m))z, \quad z(1\otimes \check{y})= (1\otimes \check{y})z, \quad m\in M, \quad \check{y}\in L^\infty(\check{\G}).
\end{align}
On the other hand, we have $z\in {}_M\mathscr{L}^\G(S^\G(\mathcal{G}))$ if and only if
    \begin{align}\label{ha}
        &z(\pi_{\mathcal{G}}\otimes \id)(\alpha(m)) = (\pi_{\mathcal{G}}\otimes \id)(\alpha(m))z, \quad  z V_{23}(\omega) = V_{23}(\omega)z, \quad m\in M, \quad \omega \in L^1(\G).
    \end{align}
    Since $[V(\omega): \omega \in L^1(\G)]^{\sigma\text{-weakly}}= L^\infty(\check{\G})$, \eqref{hah} and \eqref{ha} are equivalent. 
\end{proof}

\subsection{Generators} As in \cite{Rie74}, the following concept will be crucial in this paper:

\begin{Def}
    A generator for $\operatorname{Corr}^\G(M,N)$ consists of a $\G$-$M$-$N$-correspondence $\mathcal{H}$ such that for every $\mathcal{K}\in \Corr^\G(M,N)$, there exists a set $I$ and an intertwining isometry $\mathcal{K}\to \bigoplus_{i\in I} \mathcal{H}$. In other words, every $\mathcal{K}\in \Corr^\G(M,N)$ is contained in an amplification of $\mathcal{H}$. 
\end{Def}
We give some remarks:
\begin{itemize}
    \item Since every $\G$-$M$-$N$ correspondence decomposes as a direct sum of cyclic $\G$-$M$-$N$-correspondences \cite[Proposition 2.11]{DCDR24}, there exists a $\G$-$M$-$N$-correspondence $\mathcal{H}$ that contains (an isomorphic copy of) every cyclic $\G$-$M$-$N$-correspondence. Then $\mathcal{H}$ is a generator for $\operatorname{Corr}^\G(M,N)$, so generators for $\Corr^\G(M,N)$ always exist. 
    \item The generators of $\Rep(M)= \Corr(M, \C)$ are very easy to describe. These are exactly the faithful normal unital $*$-representations of $M$ on Hilbert spaces \cite[Proposition 1.3]{Rie74}. In particular, $(L^2(M), \pi_M)\in \Rep(M)$ is a generator for $\Rep(M)$.
    \item Let $\G$ be a compact quantum group and consider $(L^2(\G),V)\in \Rep(\G)$. 
By Peter-Weyl theory, every unitary representation decomposes as a direct sum of irreducible representations and every irreducible representation is contained in $(L^2(\G),V)$. Therefore, $(L^2(\G),V)$ is a generator for $\Rep(\G)$.
\end{itemize}

\subsection{$W^*$-categories} \cite{GLR85} A \emph{$W^*$-category} $\mathscr{C}$ is a $C^*$-category such that for every two objects $X,Y\in \mathscr{C}$, the Banach space $\mathscr{C}(X,Y)$ of morphisms $X\to Y$ has a predual. In particular, the morphism spaces $\mathscr{C}(X)=\mathscr{C}(X,X)$ are $W^*$-algebras. A functor $F: \mathscr{C}\to \mathscr{D}$ between $C^*$-categories is called \emph{$*$-functor} if $F(f^*) = F(f)^*$ for all $X,Y \in \mathscr{C}$ and all $f\in \mathscr{C}(X,Y)$. A $*$-functor $F: \mathscr{C}\to \mathscr{D}$ between $W^*$-categories is called \emph{normal} if for every object $X\in \mathscr{C}$, the unital $*$-homomorphism
$$\mathscr{C}(X)\to \mathscr{D}(F(X)): f \mapsto F(f)$$
is normal. Two $W^*$-categories $\mathscr{C}, \mathscr{D}$ are called equivalent (as $W^*$-categories) if there is a fully faithful $*$-functor $F: \mathscr{C}\to \mathscr{D}$ that is essentially surjective. Such a functor is automatically normal. Equivalently, there exist normal $*$-functors $F: \mathscr{C}\to \mathscr{D}$ and $G: \mathscr{D}\to \mathscr{C}$ such that $F\circ G \cong \id_{\mathscr{D}}$ and $G\circ F \cong \id_{\mathscr{C}}.$ All $W^*$-categories in this paper will be very concrete. A leading role will be played by the $W^*$-category $\Corr^\G(M,N)$ where $\G$ is a locally compact quantum group and $M,N$ are $\G$-$W^*$-algebras.

Let $\mathscr{C}$ be a $W^*$-category and let $\mathscr{D}$ be a $C^*$-tensor category \cite{NT14}. If $\mathscr{C}$ also obtains the structure of a (right) $\mathscr{D}$-module $C^*$-category \cite{DCY13}, then we call $\mathscr{C}$ a (right) $\mathscr{D}$-module $W^*$-category. 


\section{Functors associated to pairs of dynamical systems}

Fix throughout this entire section a locally compact quantum group $\G$ and $\G$-$W^*$-algebras $(M, \alpha), (N, \beta)$. 

Clearly, the category $\Rep(\G)$ is a $C^*$-tensor category. Given $(\mathcal{H},\pi_{\mathcal{H}}, U_{\mathcal{H}})\in \Rep^\G(M)$ and $(\mathcal{K}, U_{\mathcal{K}})\in \Rep(\G)$, it is easily verified that
$$\mathcal{H}\otimes \mathcal{K}=(\mathcal{H}\otimes \mathcal{K}, \pi_{\mathcal{H}\otimes \mathcal{K}}: M\to B(\mathcal{H}\otimes \mathcal{K}): m \mapsto \pi_{\mathcal{H}}(m)\otimes 1, U_{\mathcal{H},13}U_{\mathcal{K},23}) \in \Rep^\G(M).$$
Therefore, we have a canonical $*$-functor
$$\Rep^\G(M)\times \Rep(\G)\to \Rep^\G(M)$$
through which $\Rep^\G(M)$ becomes a right $\Rep(\G)$-module $W^*$-category. We will treat this module category as being strict, i.e. we will identify
$$\mathcal{H}\otimes \mathcal{K}\otimes \mathcal{L}:=\mathcal{H}\otimes (\mathcal{K}\otimes \mathcal{L})= (\mathcal{H}\otimes \mathcal{K})\otimes \mathcal{L}, \quad \mathcal{H}\in \Rep^\G(M), \quad \mathcal{K}, \mathcal{L}\in \Rep(\G)$$
without further mention. 

For future use, we also note that if $\mathcal{G}\in \Rep^\G(M)$, then 
$U_{\mathcal{G}} \in {}_M\mathscr{L}^\G(\mathcal{G}\otimes L^2(\G), S^\G(\mathcal{G}))$
defines a unitary isomorphism in $\Rep^\G(M)$.

We will be interested in $*$-functors compatible with the $\Rep(\G)$-module structure. The following definitions agree with the ones in \cite[Section 2.3]{DCY13} (up to a left versus right switch in conventions). 

\begin{Def}
    A $\Rep(\G)$-module $*$-functor $\Rep^\G(N)\to \Rep^\G(M)$ consists of a $*$-functor $F: \Rep^\G(N)\to \Rep^\G(M)$ and a collection of unitaries 
    $$T_{\mathcal{H}, \mathcal{K}}: F(\mathcal{H}\otimes \mathcal{K})\to F(\mathcal{H})\otimes \mathcal{K}, \quad \mathcal{H}\in \Rep^\G(N), \quad \mathcal{K}\in \Rep(\G)$$
    in $\Rep^\G(M)$,  natural in  $\mathcal{H}\in \Rep^\G(N)$ and natural in $\mathcal{K}\in \Rep(\G)$ such that the following properties hold:
    \begin{itemize}
       
        \item The diagram  \begin{equation}\label{commutativity}
\begin{tikzcd}
                                                                                                                            & F(\mathcal{H}\otimes \mathcal{K}\otimes \mathcal{L}) \arrow[ld, "{T_{\mathcal{H\otimes \mathcal{K}, \mathcal{L}}}}"'] \arrow[rd, "{T_{\mathcal{H}, \mathcal{K}\otimes \mathcal{L}}}"] &                                                      \\
F(\mathcal{H}\otimes \mathcal{K})\otimes \mathcal{L} \arrow[rr, "{T_{\mathcal{H}, \mathcal{K}}\otimes \id_{\mathcal{L}}}"'] &                                                                                                                                                                                       & F(\mathcal{H})\otimes \mathcal{K}\otimes \mathcal{L}
\end{tikzcd}\end{equation}
commutes for all $\mathcal{H}\in \Rep^\G(N)$ and $\mathcal{K}, \mathcal{L}\in \Rep(\G)$ (where, as usual, we view $\mathcal{K}\otimes \mathcal{L}$ with the tensor product $\G$-representation $U_{\mathcal{K},13}U_{\mathcal{L},23}.$)
\item The diagram \begin{equation}\label{commutativity2} 
\begin{tikzcd}
F(\mathcal{H}\otimes \mathbb{C}) \arrow[rr, "{T_{\mathcal{H}, \mathbb{C}}}"] &                                                                                   & F(\mathcal{H})\otimes \mathbb{C} \\
                                                                             & F(\mathcal{H}) \arrow[lu, "F(i_{\mathcal{H}})"] \arrow[ru, "i_{F(\mathcal{H})}"'] &                                 
\end{tikzcd}\end{equation}
commutes for all $\mathcal{H}\in \Rep^\G(N)$, where $i_{\mathcal{H}}: \mathcal{H}\to \mathcal{H}\otimes \mathbb{C}$ is the natural isomorphism.
    \end{itemize}
\end{Def}

\begin{Def}
    The $\Rep(\G)$-module $W^*$-categories $\Rep^\G(M)$ and $\Rep^\G(N)$ are called equivalent (as $\Rep(\G)$-module $W^*$-categories) if there exists a $\Rep(\G)$-module $*$-functor $F: \Rep^\G(N)\to \Rep^\G(M)$ such that $F$ is fully faithful and essentially surjective.
\end{Def}

Also the following category will be of interest to us:

\begin{Def}
    We will write $\operatorname{Fun}_{\Rep(\G)}(\Rep^\G(N), \Rep^\G(M))$ for the category defined by the following data:
\begin{enumerate}
    \item Objects consist of normal $\Rep(\G)$-module $*$-functors $F: \Rep^\G(N)\to \Rep^\G(M)$.
    \item Morphisms between two such functors $F,G: \Rep^\G(N)\to \Rep^\G(M)$ consist of natural transformations $s: F \implies G$ such that the diagram
\begin{equation}\label{tensor}
    \begin{tikzcd}
F(\mathcal{H}\otimes \mathcal{K}) \arrow[rr, "s_{\mathcal{H}\otimes \mathcal{K}}"] \arrow[d, "{T_{\mathcal{H}, \mathcal{K}}}"'] &  & G(\mathcal{H}\otimes \mathcal{K}) \arrow[d, "{T_{\mathcal{\mathcal{H, \mathcal{K}}}}}"] \\
F(\mathcal{H})\otimes \mathcal{K} \arrow[rr, "s_{\mathcal{H}}\otimes \id_{\mathcal{K}}"']                                       &  & G(\mathcal{H})\otimes \mathcal{K}                                                      
\end{tikzcd}
\end{equation}
commutes for all $\mathcal{H}\in \Rep^\G(N)$ and all $\mathcal{K}\in \Rep(\G)$. We denote this space of morphisms by $\operatorname{Nat}_{\Rep(\G)}(F,G).$ A morphism $s: F \implies G$ is called unitary isomorphism if $s_{\mathcal{H}}$ is a unitary for all $\mathcal{H}\in \Rep^\G(N)$.
\item If $F,G,H: \Rep^\G(N)\to \Rep^\G(M)$ are objects and $s: F\implies G$ and $t: G \implies H$ morphisms, the composition  $ts: F \implies H$ is defined by the object-wise composition:
$$(ts)_\mathcal{H}:= t_{\mathcal{H}}\circ s_{\mathcal{H}}, \quad \mathcal{H}\in \Rep^\G(N).$$
\end{enumerate}
\end{Def}

We will now discuss some standard examples:

\begin{Exa}\label{discussion}
    Given $\mathcal{G}\in \Corr^\G(M, N)$, consider the functor
    $F_{\mathcal{G}}: \Rep^\G(N)\to \Rep^\G(M)$
    defined as follows:
    \begin{itemize}
        \item If $\mathcal{H}\in \Rep^\G(N)$, we put $F_\mathcal{G}(\mathcal{H}):= \mathcal{G}\boxtimes_N \mathcal{H}\in \Corr^\G(M, \mathbb{C})= \Rep^\G(M).$
        \item If $\mathcal{H}, \mathcal{H}'\in \Rep^\G(N)$ and $x\in {}_N\mathscr{L}^\G(\mathcal{H}, \mathcal{H}')$, we define
        $F_{\mathcal{G}}(x)\in {}_M\mathscr{L}^\G(\mathcal{G}\boxtimes_N \mathcal{H}, \mathcal{G}\boxtimes_N \mathcal{H}')$ by
        $$F_{\mathcal{G}}(x)(y\otimes_N \xi) = y\otimes_N x\xi, \quad y \in \mathscr{L}_N(L^2(N), \mathcal{G}), \quad \xi \in \mathcal{H}.$$
    \end{itemize}
 Given $\mathcal{H}\in \Rep^\G(N), \mathcal{K}\in \Rep(\G)$, define the unitaries
    $T_{\mathcal{H}, \mathcal{K}}\in {}_M\mathscr{L}^{\G}(\mathcal{G}\boxtimes_N (\mathcal{H}\otimes \mathcal{K}),(\mathcal{G}\boxtimes_N \mathcal{H})\otimes \mathcal{K})$ by
    $$T_{\mathcal{H}, \mathcal{K}}(y\otimes_N (\xi\otimes \eta)) = (y\otimes_N \xi)\otimes \eta, \quad y \in \mathscr{L}_N(L^2(N), \mathcal{G}), \quad \xi \in \mathcal{H}, \eta \in \mathcal{K}.$$
    It is easily verified that the above data defines a normal $\Rep(\G)$-module $*$-functor $\Rep^\G(N)\to \Rep^\G(M)$.

    Given $\mathcal{H}, \mathcal{G}\in \Corr^\G(M,N)$, consider an intertwiner $x\in {}_M\mathscr{L}^\G_N(\mathcal{H}, \mathcal{G})$. If $\mathcal{L}\in \Rep^\G(N)$, consider the morphism
$$P(x)_\mathcal{L}:= x \boxtimes_N \id \in {}_M\mathscr{L}^\G(\mathcal{H}\boxtimes_N \mathcal{L},\mathcal{G}\boxtimes_N \mathcal{L}), \quad P(x)_{\mathcal{L}}(y\otimes_N \xi)= xy\otimes_N \xi, \quad y \in \mathscr{L}_N(L^2(N), \mathcal{H}), \quad \xi \in \mathcal{L}.$$
Then we see that
$$P(x):= (P(x)_{\mathcal{L}})_{\mathcal{L}\in \Rep^\G(N)}\in \operatorname{Nat}_{\Rep(\G)}(F_{\mathcal{H}}, F_{\mathcal{G}}).$$

\end{Exa}


From the discussion in Example \ref{discussion}, we obtain a canonical functor
$$P: \operatorname{Corr}^\G(M,N)\to \operatorname{Fun}_{\Rep(\G)}(\Rep^\G(N), \Rep^\G(M)).$$
The main goal of this section is to construct a canonical functor
$$Q: \operatorname{Fun}_{\Rep(\G)}(\Rep^\G(N), \Rep^\G(M))\to  \operatorname{Corr}^\G(M,N).$$
The following result explains how to define this functor on objects:
\begin{Theorem}\label{mr} Let $\G$ be a locally compact quantum group.
    Suppose that $F: \Rep^\G(N)\to \Rep^\G(M)$ is a normal $\Rep(\G)$-module $*$-functor. Consider $(\mathcal{G}, \pi_{\mathcal{G}}, U_{\mathcal{G}}):=  F(L^2(N))\in \Rep^\G(M)$. The Hilbert space $\mathcal{G}$ admits a canonical anti-$*$-representation
    $\rho_{\mathcal{G}}: N \to B(\mathcal{G})$
    such that $(\mathcal{G}, \pi_{\mathcal{G}}, \rho_{\mathcal{G}}, U_{\mathcal{G}})\in \Corr^\G(M,N).$
\end{Theorem}

Before diving into the proof, we isolate the following necessary lemma:

\begin{Lem}\label{intermediary} Let $\G$ be a locally compact quantum group.
Let $(\mathcal{H}, \pi, U)\in \Rep^\G(M)$ and let $\check{\rho}: N\rtimes_\beta \G\to B(\mathcal{H}\otimes L^2(\G))$ be a unital normal $*$-anti-representation such that 
$$(\mathcal{H}\otimes L^2(\G), \pi^\rtimes, \check{\rho}, \check{V}_{23})\in \operatorname{Corr}^{\check{\G}}(M\rtimes_\alpha \G, N\rtimes_\beta \G)$$ and such that $\check{\rho}(1\otimes \check{y})= U(1\otimes \rho_{L^\infty(\check{\G})}(\check{y}))U^*$ for all $\check{y}\in L^\infty(\check{\G}).$ Then there is a unique anti-$*$-representation $\rho: N \to B(\mathcal{H})$ such that $(\mathcal{H}, \pi, \rho, U) \in \Corr^\G(M,N)$ and $\rho^\rtimes = \check{\rho}$. 
\end{Lem}
\begin{proof}
   Given $n\in N$, we have 
    $\check{\rho}(\beta(n))\otimes 1= \check{V}_{23}^*(\check{\rho}(\beta(n))\otimes 1)\check{V}_{23}$,
    from which it follows that $\check{\rho}(\beta(n))\in B(\mathcal{H})\ovot L^\infty(\G)$. On the other hand, from the fact that the images of $\pi^\rtimes$ and $\check{\rho}$ commute, we also have
    $$\check{\rho}(\beta(n))(1\otimes \check{y}) =  (1\otimes \check{y})\check{\rho}(\beta(n)), \quad \check{y}\in L^\infty(\check{\G}).$$
    
    Therefore, also $\check{\rho}(\beta(n))\in B(\mathcal{H})\ovot L^\infty(\hat{\G})$. We conclude that $\check{\rho}(\beta(n))\in B(\mathcal{H})\ovot \mathbb{C}$,
    so there is a unique normal anti-$*$-representation $\rho: N\to B(\mathcal{H})$
    such that
    $$\check{\rho}(\beta(n)) = \rho(n)\otimes 1, \quad n \in N.$$
    Given $m\in M$ and $n\in N$, the commutation relation
    $$(\rho(n)\otimes 1)(\pi\otimes \id)(\alpha(m)) = (\pi\otimes \id)(\alpha(m))(\rho(n)\otimes 1)$$
    implies that 
    $$\rho(n)\pi((\id \otimes \omega)\alpha(m)) = \pi((\id \otimes \omega)\alpha(m)) \rho(n), \quad \omega \in L^1(\G).$$
    By weak Podleś density \cite[Corollary 2.7]{KS15}, 
    $$\rho(n)\pi(m) = \pi(m)\rho(n), \quad m\in M, n \in N.$$
     It remains to show that $\rho$ is compatible with $U$, i.e.
    $$(\rho\otimes R)\beta(n) = U^*(\rho(n)\otimes 1)U, \quad n \in N.$$
Note first that
$(\check{J}\otimes \check{J})V^*(\check{J}\otimes \check{J})= W_{21}$
so that $(\check{\rho}\otimes R)(V_{23})= U_{12}W_{32}U_{12}^* = W_{32}U_{13}.$
Therefore, we find
\begin{align*}
    U_{13}^*(\rho(n)\otimes 1\otimes 1)U_{13}&= U_{13}^*W_{32}^*(\rho(n)\otimes 1\otimes 1)W_{32}U_{13}= (\check{\rho}\otimes R)(V_{23}^*)(\check{\rho}\otimes R)(\beta(n)\otimes 1))(\check{\rho}\otimes R)(V_{23})\\
    &=(\check{\rho}\otimes R)(V_{23}(\beta(n)\otimes 1)V_{23}^*)
    = (\check{\rho}\otimes R)(\id \otimes \Delta)(\beta(n))\\
    &= (\check{\rho}\otimes R)(\beta \otimes \id)(\beta(n))= (\rho\otimes R)(\beta(n))_{13}
\end{align*}
from which it follows that
$U^*(\rho(n)\otimes 1)U= (\rho \otimes R)\beta(n)$,
as desired.
\end{proof}

\begin{proof}[Proof of Theorem \ref{mr}]  Consider the unitary isomorphism
   $$
\begin{tikzcd}
S: F(S_N^\G) \arrow[rr, "F(U_\beta^*)"] &  & F(L^2(N)\otimes L^2(\G)) \arrow[rr, "{T_{L^2(N), L^2(\G)}}"] &  & \mathcal{G}\otimes L^2(\G) \arrow[rr, "U_{\mathcal{G}}"] &  & S^\G(\mathcal{G})
\end{tikzcd}$$
in $\Rep^\G(M)$. By Lemma \ref{easy}, we have
$$\rho_{N\rtimes_\beta \G}(N\rtimes_\beta \G)= \pi_{N\rtimes_\beta \G}(N\rtimes_\beta \G)'= {}_N\mathscr{L}^\G(S_N^\G).$$
Therefore, it makes sense to consider the following data on the Hilbert space $\mathcal{G}\otimes L^2(\G)$:
\begin{itemize}
    \item The normal $*$-representation $\pi_{\mathcal{G}}^\rtimes: M\rtimes_\alpha \G \to B(\mathcal{G}\otimes L^2(\G)): z \mapsto (\pi_{\mathcal{G}}\otimes \id)(z).$
    \item The normal $*$-anti-representation $\check{\rho}: N \rtimes_\beta \G \to B(\mathcal{G}\otimes L^2(\G)): z \mapsto S F(\rho_{N\rtimes \G}(z))S^*$.
\end{itemize}
Then we have
\begin{align}\label{refer}
    \check{\rho}(N\rtimes_\beta \G)= S F({}_N\mathscr{L}^\G(S_N^\G))S^* \subseteq S [{}_M \mathscr{L}^\G(F(S_N^\G))]S^*= {}_M \mathscr{L}^\G(S^\G(\mathcal{G})) = \pi_{\mathcal{G}}^\rtimes(M\rtimes_\alpha \G)'
\end{align}
so that $(\mathcal{G}\otimes L^2(\G), \pi_{\mathcal{G}}^\rtimes, \check{\rho})\in \Corr(M\rtimes_\alpha \G, N\rtimes_\beta \G).$

Note that for $\check{y}\in L^\infty(\check{\G})$, we have that $\rho_{L^\infty(\check{\G})}(\check{y}) \in \mathscr{L}^\G(L^2(\G))$, so by naturality 
\begin{align*}
    \check{\rho}(1\otimes \check{y}) &= SF(U_\beta(1\otimes \rho_{L^\infty(\check{\G})}(\check{y}))U_\beta^*)S^*\\
    &= U_{\mathcal{G}}T_{L^2(N), L^2(\G)}F(1\otimes \rho_{L^\infty(\check{\G})}(\check{y}))T_{L^2(N), L^2(\G)}^* U_{\mathcal{G}}^*= U_{\mathcal{G}}(1\otimes \rho_{L^\infty(\check{\G})}(\check{y}))U_{\mathcal{G}}^*
\end{align*}
and for $n\in N$,
\begin{align*}
    \check{\rho}(\beta(n)) = SF(\rho_N(n)\otimes 1)S^* = U_{\mathcal{G}} T_{L^2(N), L^2(\G)} F((\rho_N\otimes R)\beta(n)) T_{L^2(N), L^2(\G)}^* U_{\mathcal{G}}^*.
\end{align*}

We now show that $(\mathcal{G}\otimes L^2(\G), \pi^\rtimes_{\mathcal{G}}, \check{\rho}, \check{V}_{23})\in \Corr^{\check{\G}}(M\rtimes_\alpha \G, N\rtimes_\beta \G)$, i.e. we must show that
$$(\check{\rho}\otimes \check{R})(\id \otimes \check{\Delta}_r)(z) =\check{V}_{23}^*(\check{\rho}(z)\otimes 1)\check{V}_{23}, \quad z \in N\rtimes_\beta \G.$$
It suffices to show this for for $z= 1\otimes \check{y}$ with $\check{y}\in L^\infty(\check{\G})$ and for $z= \beta(n)$ with $n\in N$. In the first case,
  let $\check{y} \in L^\infty(\check{\G})$. Then we have
    \begin{align*}
    (\rho_{L^\infty(\check{\G})}\otimes \check{R})\check{\Delta}(\check{y}) = (\check{J}\otimes J)\check{\Delta}(\check{y}^*)(\check{J}\otimes J) = \check{V}^*(\rho_{L^\infty(\check{\G})}(\check{y})\otimes 1)\check{V}.
    \end{align*}
    Therefore, if $\check{\mu} \in L^1(\check{\G})$, we find 
    \begin{align*}
        (\id \otimes \id \otimes \check{\mu})(\check{\rho}\otimes \check{R})(1\otimes \check{\Delta}(\check{y}))&= \check{\rho}(1\otimes (\id \otimes \check{\mu}\check{R})\check{\Delta}(\check{y}))\\
        &= U_{\mathcal{G}}(1\otimes \rho_{L^\infty(\check{\G})}((\id \otimes \check{\mu}\check{R})\check{\Delta}(\check{y})))U_{\mathcal{G}}^*\\
        &=U_{\mathcal{G}}(1\otimes (\id \otimes \check{\mu})(\check{V}^*(\rho_{L^\infty(\check{\G})}(\check{y})\otimes 1)\check{V}))U_{\mathcal{G}}^*\\
        &= (\id \otimes \id \otimes \check{\mu})(\check{V}_{23}^* U_{\mathcal{G},12}(1\otimes \rho_{L^\infty(\check{\G})}(\check{y})\otimes 1 )U_{\mathcal{G},12}^*\check{V}_{23})\\
        &= (\id \otimes \id \otimes \check{\mu})(\check{V}_{23}^* (\check{\rho}(1\otimes \check{y})\otimes 1)\check{V}_{23})
    \end{align*}
    from which we deduce that 
    $$(\check{\rho}\otimes \check{R})(1\otimes \check{\Delta}(\check{y})) = \check{V}_{23}^*(\check{\rho}(1\otimes \check{y})\otimes 1)\check{V}_{23}.$$

In the second case, we must show that 
$$\check{\rho}(\beta(n))\otimes 1=\check{V}_{23}^*(\check{\rho}(\beta(n))\otimes 1)\check{V}_{23}$$
or equivalently
\begin{equation}\label{equ}
    T_{L^2(N), L^2(\G)} F(U_\beta^*(\rho_N(n)\otimes 1)U_\beta)T_{L^2(N), L^2(\G)}^*\in B(\mathcal{G})\ovot L^\infty(\G).
\end{equation}

Let us now simplify some notation. We put 
\begin{align*}
    T_V&:= T_{L^2(N), L^2(\G)}= T_{L^2(N), (L^2(\G), V)}\\
    T_{V_{13}}&:= T_{L^2(N), (L^2(\G)\otimes L^2(\G), V_{13})}\\
    T_{V_{13}W_{32}}&:= T_{L^2(N), (L^2(\G)\otimes L^2(\G), V_{13}W_{32})}\\
    T_{\mathbb{I}}^{\otimes}&:= T_{L^2(N)\otimes L^2(\G), (L^2(\G), \mathbb{I})}\\
    T_{W_{21}}^\otimes&:= T_{L^2(N)\otimes L^2(\G), (L^2(\G), W_{21})}.
\end{align*}

Note that by \eqref{commutativity}, we have that
$$T_{V_{13}}\circ T_{\mathbb{I}}^{\otimes*} = T_V\otimes \id_{L^2(\G)}= T_{V_{13}W_{32}}\circ T_{W_{21}}^{\otimes *}$$

Put $U:= (u_\G\otimes 1)W^*(u_\G\otimes 1)\in L^\infty(\G)' \ovot L^\infty(\hat{\G})$. Then an easy computation shows that 
$V_{13}= U_{12}V_{13}W_{32}U_{12}^*.$
In other words, $U$ defines a unitary isomorphism
$$U: (L^2(\G)\otimes L^2(\G), V_{13}W_{32})\to (L^2(\G)\otimes L^2(\G), V_{13})$$
in the category $\Rep(\G)$. Therefore, by naturality of the unitaries $T_{\mathcal{H}, \mathcal{K}}$ in the second variable, we deduce that
\begin{equation}
    (\id_{\mathcal{G}}\otimes U)\circ T_{V_{13}W_{32}}= T_{V_{13}}\circ F(\id_{L^2(N)}\otimes U).
\end{equation}

Consequently, since $U_\beta^*(\rho_N(n)\otimes 1)U_\beta \in {}_N\mathscr{L}^\G(L^2(N)\otimes L^2(\G))$, we find  
\begin{align*}
&(\id_{\mathcal{G}}\otimes U)(T_V\otimes \id_{L^2(\G)})(F(U_\beta^*(\rho_N(n)\otimes 1)U_\beta)\otimes \id_{L^2(\G)})(T_V^*\otimes \id_{L^2(\G)})(\id_{\mathcal{G}}\otimes U^*)\\
&= (\id_{\mathcal{G}}\otimes U) T_{V_{13}W_{32}} T_{W_{21}}^{\otimes *} (F(U_\beta^*(\rho_N(n)\otimes 1)U_\beta)\otimes \id_{L^2(\G)}) T_{W_{21}}^{\otimes} T_{V_{13}W_{32}}^* (\id_{\mathcal{G}}\otimes U^*)\\
&= (\id_{\mathcal{G}}\otimes U)T_{V_{13}W_{32}} F(U_\beta^*(\rho_N(n)\otimes 1)U_\beta \otimes \id_{L^2(\G)}) T_{V_{13}W_{32}}^*(\id_{\mathcal{G}}\otimes U^*)\\
&= T_{V_{13}}F(\id_{L^2(N)}\otimes U)F(U_{\beta,12}^*(\rho_N(n)\otimes 1\otimes 1)U_{\beta,12}) F(\id_{L^2(N)}\otimes U^*)T_{V_{13}}^*\\
&= T_{V_{13}} F(U_{\beta,12}^*(\rho_N(n)\otimes 1 \otimes 1)U_{\beta,12})T_{V_{13}}^*\\
&= T_{V_{13}}T_{\mathbb{I}}^{\otimes*}(F(U_\beta^*(\rho_N(n)\otimes 1)U_\beta)\otimes \id_{L^2(\G)}) T_{\mathbb{I}}^{\otimes} T_{V_{13}}^*\\
&= (T_V\otimes \id_{L^2(\G)})(F(U_\beta^*(\rho_N(n)\otimes 1)U_\beta)\otimes \id_{L^2(\G)})(T_V^*\otimes \id_{L^2(\G)}).
\end{align*}

From the fact that
$[(\id \otimes \omega)(U): \omega \in B(L^2(\G))_*]^{\sigma\text{-weak}}= L^\infty(\G)'$
we deduce that \eqref{equ} holds. Consequently, $(\mathcal{G}\otimes L^2(\G), \pi^\rtimes_{\mathcal{G}}, \check{\rho}, \check{V}_{23})\in \Corr^{\check{\G}}(M\rtimes_\alpha \G, N\rtimes_\beta \G)$. Therefore, the result follows from Lemma \ref{intermediary}.
\end{proof}

Next, we explain how to define the functor 
$$Q: \operatorname{Fun}_{\Rep(\G)}(\Rep^\G(N), \Rep^\G(M))\to \Corr^\G(M,N)$$
on morphisms.
\begin{Prop}\label{functors}
    Let $\G$ be a locally compact quantum group, $F,G\in \operatorname{Fun}_{\Rep(\G)}(\Rep^\G(N), \Rep^\G(M))$ and $s\in \operatorname{Nat}_{\Rep(\G)}(F,G)$. 
    Then $s_{L^2(N)} \in {}_M\mathscr{L}_N^\G(F(L^2(N)), G(L^2(N))).$ 
\end{Prop}
  
\begin{proof}
  Recalling that $(\rho_N\otimes R)\beta(n) \in {}_N\mathscr{L}^\G(L^2(N)\otimes L^2(\G))$ for every $n\in N$, we calculate
    \begin{align*} s_{L^2(N)} \rho_{F(L^2(N))}(n)\otimes 1&=
        (s_{L^2(N)}\otimes \id)\rho_{F(L^2(N))}^\rtimes(\beta(n))\\
        &= (s_{L^2(N)}\otimes \id)U_{F(L^2(N))} T_{L^2(N), L^2(\G)} F((\rho_N\otimes R)\beta(n))T_{L^2(N), L^2(\G)}^* U_{F(L^2(N))}^*\\
        &= U_{G(L^2(N))}(s_{L^2(N)}\otimes \id) T_{L^2(N), L^2(\G)} F((\rho_N\otimes R)\beta(n)) T_{L^2(N), L^2(\G)}^* U_{F(L^2(N))}^*\\
        &= U_{G(L^2(N))} T_{L^2(N), L^2(\G)} s_{L^2(N)\otimes L^2(\G)} F((\rho_N\otimes R)\beta(n))T_{L^2(N), L^2(\G)}^* U_{F(L^2(N))}^*\\
        &= U_{G(L^2(N))} T_{L^2(N), L^2(\G)} G((\rho_N\otimes R)\beta(n)) s_{L^2(N)\otimes L^2(\G)} T_{L^2(N), L^2(\G)}^* U_{F(L^2(N))}^*\\
        &= U_{G(L^2(N))} T_{L^2(N), L^2(\G)} G((\rho_N\otimes R)\beta(n)) T_{L^2(N), L^2(\G)}^* (s_{L^2(N)}\otimes \id) U_{F(L^2(N))}^*\\
        &= U_{G(L^2(N))} T_{L^2(N), L^2(\G)} G((\rho_N\otimes R)\beta(n)) T_{L^2(N), L^2(\G)}^* U_{G(L^2(N))}^*(s_{L^2(N)}\otimes \id) \\
        &= \rho_{G(L^2(N))}^{\rtimes}(\beta(n))(s_{L^2(N)}\otimes \id) = \rho_{G(L^2(N))}(n) s_{L^2(N)}\otimes 1.
    \end{align*}
  The desired conclusion follows.
\end{proof}

Combining Theorem \ref{mr} and Proposition \ref{functors}, we obtain a well-defined functor
$$Q: \operatorname{Fun}_{\Rep(\G)}(\Rep^\G(N), \Rep^\G(M))\to \operatorname{Corr}^\G(M,N).$$
It is easily observed that $Q\circ P\cong \id_{\Corr^\G(M,N)}$. Indeed, if $\mathcal{G}\in \Corr^\G(M,N)$, we have a unitary isomorphism
$$QP(\mathcal{G})= \mathcal{G}\boxtimes_N L^2(N)\to \mathcal{G}: y \otimes_N \xi \mapsto y\xi, \quad y \in \mathscr{L}_N(L^2(N), \mathcal{G}), \quad \xi \in L^2(N)$$
of $\G$-$M$-$N$-correspondences. These isomorphisms are easily checked to be natural in $\mathcal{G}\in \Corr^\G(M,N)$. 

Next, we give a categorical characterisation of equivariant Morita equivalence of dynamical von Neumann algebras. First, we need two lemmas:

\begin{Lem}\label{inclusion}
    If $(\mathcal{H}, \kappa)\in \Rep(M\rtimes_\alpha \G)$, there exists a unique $(\mathcal{H}, \pi_\kappa, U_\kappa)\in \Rep^\G(M)$ such that 
$$\kappa(\alpha(m)) = \pi_\kappa(m), \quad \kappa(1\otimes V(\omega)) = U_\kappa(\omega), \quad m\in M, \quad \omega \in L^1(\G).$$
In this way, we obtain a canonical fully faithful normal $*$-functor $\Rep(M\rtimes_\alpha \G)\to \Rep^\G(M)$ that is the identity on morphisms.
\end{Lem}
\begin{proof}
Given a normal unital $*$-homomorphism
    $\kappa: M\rtimes_\alpha \G \to B(\mathcal{H})$, define $\pi_\kappa:= \kappa \circ \alpha: M\to B(\mathcal{H})$. There is a unique unitary $\G$-representation $U_\kappa \in B(\mathcal{H})\ovot L^\infty(\G)$ such that
    $$\kappa(1\otimes V(\omega)) = U_\kappa(\omega), \quad \omega \in L^1(\G).$$  
    In particular, $(\kappa\otimes \id)(V_{23})= U$. We then calculate for $m\in M$ and $\omega \in L^1(\G)$ that
    \begin{align*}
        (\id \otimes \omega)(U_\kappa(\pi_\kappa(m)\otimes 1))&= U_\kappa(\omega)\pi_\kappa(m)= \kappa((1\otimes V(\omega))\alpha(m))= \kappa((\id \otimes \id \otimes \omega)(V_{23}(\alpha(m)\otimes 1)))\\
        &= \kappa((\id \otimes \id \otimes \omega)((\alpha\otimes \id)\alpha(m) V_{23}))\\
        &= (\id \otimes \omega)(\kappa\otimes \id)((\alpha\otimes \id)\alpha(m)V_{23})= (\id \otimes \omega) ((\pi_\kappa\otimes \id)(\alpha(m)) U_\kappa).
    \end{align*}
From this, we conclude that $(\mathcal{H}, \pi_\kappa, U_\kappa)\in \Rep^\G(M)$.    
\end{proof}

We will use Lemma \ref{inclusion} to view $\Rep(M\rtimes_\alpha \G)\subseteq \Rep^\G(M).$ If $\mathcal{H}\in \Rep^\G(M)$ and $L^2(\G)\in \Rep(\G)$, then 
$\mathcal{H}\otimes L^2(\G)\in \Rep^\G(M)$ actually belongs to $\Rep(M\rtimes_\alpha \G)$. Indeed, it corresponds to the $*$-representation
$$M\rtimes_\alpha \G \to B(\mathcal{H}\otimes L^2(\G)): z\mapsto U_\mathcal{H}^*(\pi_\mathcal{H}\otimes \id)(z) U_\mathcal{H}.$$
In particular, $\mathcal{H}\otimes L^2(\G)$ is a generator for $\Rep(M\rtimes_\alpha \G)$ if and only if $\pi_{\mathcal{H}}: M \to B(\mathcal{H})$ is faithful.
 
\begin{Lem}\label{induced}
    A functor $F\in \operatorname{Fun}_{\Rep(\G)}(\Rep^\G(N), \Rep^\G(M))$ restricts to a normal $*$-functor
    $$F: \Rep(N\rtimes_\beta \G)\to \Rep(M\rtimes_\alpha \G).$$ 
\end{Lem}
\begin{proof}
    Let $\mathcal{H} \in \Rep(N\rtimes_\beta \G)$. Since $L^2(N)\otimes L^2(\G)$ is a generator for $\Rep(N\rtimes_\beta \G)$, there is an isometry
    $$\iota \in {}_{N\rtimes_\beta \G}\mathscr{L}\left(\mathcal{H}, \bigoplus_{i\in I} (L^2(N)\otimes L^2(\G))\right)= {}_{N}\mathscr{L}^\G\left(\mathcal{H}, \bigoplus_{i\in I} (L^2(N)\otimes L^2(\G))\right).$$
  It induces the isometry
    $$\lambda: F(\mathcal{H})\hookrightarrow F\left(\bigoplus_{i\in I} (L^2(N)\otimes L^2(\G))\right)\cong \bigoplus_{i\in I} F(L^2(N)\otimes L^2(\G))\cong \bigoplus_{i\in I} (F(L^2(N))\otimes L^2(\G))$$
    where we made use of the fact that $F$ preserves direct sums, cfr.\ \cite[Proposition 4.9]{Rie74}.
    Then define the normal $*$-homomorphism
    $$M\rtimes_\alpha \G \to B(F(\mathcal{H})): z \mapsto \lambda^* \left(\bigoplus_{i\in I} U_{F(L^2(N))}^*(\pi_{F(L^2(N))}\otimes \id)(z)U_{F(L^2(N))}\right) \lambda$$
   which under the inclusion $\Rep(M\rtimes_\alpha \G)\subseteq \Rep^\G(M)$ corresponds to $F(\mathcal{H})\in \Rep^\G(M).$ Therefore, $F(\mathcal{H})\in \Rep(M\rtimes_\alpha \G)$.
\end{proof}


\begin{Theorem}\label{equivariant_morita}
    Let $\G$ be a locally compact quantum group. Then $(M, \alpha)\sim_\G (N, \beta)$ if and only if $\Rep^\G(M)$ and $\Rep^\G(N)$ are equivalent as $\Rep(\G)$-module $W^*$-categories.
\end{Theorem}
\begin{proof} If $\mathcal{G}\in \Corr^\G(M,N)$ defines a $\G$-$W^*$-Morita equivalence, then the functor $F_\mathcal{G}: \Rep^\G(N)\to \Rep^\G(M)$ is an equivalence of $\Rep(\G)$-module $W^*$-categories. Indeed, a quasi-inverse is given by
$$F_{\overline{\mathcal{G}}}: \Rep^\G(M)\to \Rep^\G(N)$$
as follows easily from \cite[Proposition 5.10]{DCDR24} and associativity of the Connes fusion tensor product.

Conversely, assume that $F: \Rep^\G(N)\to \Rep^\G(M)$ is a  normal $\Rep(\G)$-module $*$-functor that is fully faithful and essentially surjective. We prove that $F(L^2(N)) = (\mathcal{G}, \pi_{\mathcal{G}}, \rho_{\mathcal{G}}, U_{\mathcal{G}})\in \Corr^\G(M,N)$ defines an equivariant Morita equivalence.
    
    Since $F$ is faithful, the right action 
    $\rho_{\mathcal{G}}^\rtimes: \mathcal{G}\otimes L^2(\G)\curvearrowleft N\rtimes_\beta \G$ is faithful. On the other hand, note that since $F$ is full, the inclusion in \eqref{refer} must be an equality. Therefore, $\rho_{\mathcal{G}}^\rtimes(N\rtimes_\beta \G) = \pi_{\mathcal{G}}^\rtimes(M\rtimes_\alpha \G)'.$ We now argue that the left action $\pi_{\mathcal{G}}^\rtimes: M\rtimes_\alpha \G \curvearrowright \mathcal{G}\otimes L^2(\G)$ is faithful. By Lemma \ref{induced}, the functor $F$ restricts to a normal $*$-functor 
    $$F: \Rep(N\rtimes_\beta \G)\to \Rep(M\rtimes_\alpha \G)$$
 which is still an equivalence of $W^*$-categories\footnote{This is easily seen by making use of the fact that there exists a quasi-inverse $G: \Rep^\G(M)\to \Rep^\G(M)$ that is also a normal $\Rep(\G)$-module $*$-functor.}. Therefore, $F(L^2(N)\otimes L^2(\G))\cong \mathcal{G}\otimes L^2(\G)$ is a generator for $\Rep(M\rtimes_\alpha \G)$ and thus the action
 $\pi_{\mathcal{G}}^\rtimes: M\rtimes_\alpha \G \curvearrowright \mathcal{G}\otimes L^2(\G)$
 is faithful. Therefore, $\mathcal{G}^\rtimes$ is a $\check{\G}$-$W^*$-Morita equivalence. By Proposition \ref{crossed product}, $\mathcal{G}$ is a $\G$-$W^*$-Morita equivalence.
\end{proof}

\section{An equivariant Eilenberg-Watts theorem for compact quantum groups}

In the previous section, for $\G$ a locally compact quantum group and $M,N$ two $\G$-$W^*$-algebras, we constructed two canonical functors
\begin{align*}
    &P: \operatorname{Corr}^\G(M,N)\to \operatorname{Fun}_{\Rep(\G)}(\Rep^\G(N), \Rep^\G(M))\\
    &Q: \operatorname{Fun}_{\Rep(\G)}(\Rep^\G(N), \Rep^\G(M))\to \operatorname{Corr}^\G(M,N).
\end{align*}
In this section, we show that if $\G$ is a compact quantum group, then these functors are quasi-inverse to each other. In particular, every normal $\Rep(\G)$-module $*$-functor $\Rep^\G(N)\to \Rep^\G(M)$ arises as tensoring with a $\G$-$M$-$N$-correspondence, which can be seen as an equivariant version of the Eilenberg-Watts theorem.

The following proposition shows that for compact quantum groups, problems about the equivariant representation category $\Rep^\G(M)$ can be reduced to equivalent problems about the representation category $\Rep(M\rtimes_\alpha \G)$. This allows us to solve problems in the equivariant setting by appealing to the analogous non-equivariant results.

\begin{Prop}\label{compact Morita} Let $\G$ be a compact quantum group and  $(M, \alpha), (N, \beta)$ be $\G$-$W^*$-algebras.
\begin{enumerate}
    \item If $(\mathcal{H}, \pi, U)\in \Rep^\G(M)$, there exists a unique normal unital $*$-representation
$\pi\rtimes U: M\rtimes_\alpha \G \to B(\mathcal{H})$
such that 
$$(\pi\rtimes U)\alpha(m) = \pi(m), \quad (\pi\rtimes U)(1\otimes V(\omega))= U(\omega), \quad m\in M, \quad \omega \in L^1(\G).$$ 
 \item The assignment 
 $(\mathcal{H}, \pi, U)\mapsto (\mathcal{H}, \pi\rtimes U)$ defines a functor $\Rep^\G(M)\to \Rep(M\rtimes_\alpha \G)$ that is the identity on morphisms. It is an isomorphism of $W^*$-categories.
 \item $\Rep^\G(M)$ and $\Rep^\G(N)$ are equivalent $W^*$-categories if and only if $M\rtimes_\alpha \G$ and $N\rtimes_\beta \G$ are Morita equivalent von Neumann algebras.
 \item \begin{enumerate}
 \item $L^2(M)\otimes L^2(\G)$ is a generator for $\Rep^\G(M).$
 \item $L^2(M)\otimes (L^2(\G)\otimes L^2(\G))= (L^2(M)\otimes L^2(\G))\otimes L^2(\G)$ is a generator for $\Rep^\G(M)$.
 \end{enumerate}
 \item Let $\mathcal{H}$ be a generator for $\Rep^\G(N)$ and $F,G: \Rep^\G(N)\to \Rep^\G(M)$ two normal $*$-functors. If $s_1\in {}_M\mathscr{L}^\G(F(\mathcal{H}), G(\mathcal{H}))$ satisfies 
$$s_1 F(x) = G(x) s_1, \quad x \in {}_N\mathscr{L}^\G(\mathcal{H}),$$
then there exists a unique natural transformation $s\in \operatorname{Nat}(F,G)$ such that $s_{\mathcal{H}}= s_1$.
\end{enumerate} 
\end{Prop}
\begin{proof}
    (1) Since 
    $U^*(1\ovot L^\infty(\check{\G}))U\subseteq B(\mathcal{H})\ovot L^\infty(\check{\G})$ \cite[Lemma 2.5]{DR24},
we can well-define
$$\pi\rtimes U: M\rtimes_\alpha \G \to B(\mathcal{H}): z \mapsto (\id \otimes \check{\epsilon})(U^*(\pi\otimes \id)(z)U)$$
where $\check{\epsilon}\in \ell^1(\check{\G})$ is the counit of the dual discrete quantum group $\check{\G}$. We have 
    $(\check{\Delta}\otimes \id)(V) = V_{13}V_{23}$. Applying $\check{\epsilon}\otimes \id \otimes \id$ to this equality, we find $V= (1\otimes (\check{\epsilon}\otimes \id)(V))V$ hence $(\check{\epsilon}\otimes \id)(V) = 1$. Therefore, it follows that for $\omega \in L^1(\G)$
    \begin{align*}
        (\id \otimes \check{\epsilon})(U^*(1\otimes V(\omega))U) &=  (\id \otimes \check{\epsilon})(\id \otimes \id \otimes \omega)(U_{13}V_{23})\\
        &= (\id \otimes \omega)(\id \otimes \check{\epsilon}\otimes \id)(U_{13}V_{23})= (\id \otimes \omega)(U) = U(\omega).
    \end{align*} This shows that the normal $*$-homomorphism $\pi\rtimes U$ as in $(1)$ exists.
    
    $(2)$ This follows from (1) and Lemma \ref{inclusion}.
    
    $(3)$ We recall from \cite[Theorem 3.5.5]{Bro03} or  \cite[Theorem 7.9]{Rie74} that for two von Neumann algebras $A$ and $B$, we have that $\Rep(A)$ and $\Rep(B)$ are equivalent as $W^*$-categories if and only if $A$ and $B$ are Morita equivalent von Neumann algebras. The result then follows from $(2)$.

     (4) (a) In section $3$, it was already remarked that $L^2(M)\otimes L^2(\G)\in \Rep(M\rtimes_\alpha \G)= \Rep^\G(M)$ is a generator.

   (b) Write $\xi_\G = \Lambda(1) \in L^2(\G)$. Since $V(\xi_\G \otimes \xi) = \xi_\G \otimes \xi$ for all $\xi \in L^2(\G)$, it is easily verified that the map
   $$L^2(M)\otimes L^2(\G)\to L^2(M)\otimes L^2(\G)\otimes L^2(\G): \xi \otimes \eta \mapsto \xi\otimes \eta \otimes \xi_\G$$
   defines an isometry in ${}_M\mathscr{L}^\G(L^2(M)\otimes L^2(\G), L^2(M)\otimes L^2(\G)\otimes L^2(\G)).$ The fact that $L^2(M)\otimes L^2(\G)\otimes L^2(\G)$ defines a generator then follows from $(a)$. Alternatively, as already remarked in section 3, $L^2(M)\otimes (L^2(\G)\otimes L^2(\G))$ is a generator for $\Rep(M\rtimes_\alpha \G)$ with associated faithful $*$-representation
   $$M\rtimes_\alpha \G \to B(L^2(M)\otimes L^2(\G)\otimes L^2(\G)): z \mapsto V_{23}^* U_{\alpha,13}^*(\pi_M\otimes \id)(z)_{13}U_{\alpha,13}V_{23}.$$ 

(5) This follows immediately from \cite[Proposition 5.4]{Rie74} and $(2)$. \end{proof}

\begin{Rem}
    Compactness in the previous result is essential. For example, $(1)$ in Proposition \ref{compact Morita} is only true for compact quantum groups. Indeed, let $\G$ be a locally compact quantum group. Consider $\mathbb{C}\in \Rep^\G(M)$. Then we require the existence of a normal $*$-character $L^\infty(\check{\G})\to \mathbb{C}$, which already forces $\G$ to be compact.
\end{Rem}

We now provide an explicit example of a compact (quantum) group and two $\G$-$W^*$-algebras $M$ and $N$ such that $\Rep(M)$ and $\Rep(N)$ are equivalent as $W^*$-categories, but not equivalent as $\Rep(\G)$-module $W^*$-categories. By Theorem \ref{equivariant_morita} and Proposition \ref{compact Morita}, it is sufficient to find an example of two $\G$-$W^*$-algebras $(M, \alpha), (N, \beta)$ such that $(M, \alpha)\not\sim_\G (N, \beta)$ and $M\rtimes_\alpha \G \sim N \rtimes_\beta \G$:

\begin{Exa} 
Consider any non-trivial discrete (quantum) group $\mathbbl{\Gamma}.$ Then we see that
    $(\ell^\infty(\mathbbl{\Gamma}), \Delta)\not\sim_{\mathbbl{\Gamma}} (\ell^\infty(\mathbbl{\Gamma}), \tau)$. Indeed, taking $\mathbbl{\Gamma}$-crossed products, we would otherwise find that
    $$\mathbb{C}\sim B(\ell^2(\mathbbl{\Gamma}))\cong \ell^\infty(\mathbbl{\Gamma})\rtimes_\Delta \mathbbl{\Gamma} \sim \ell^\infty(\mathbbl{\Gamma})\rtimes_\tau \mathbbl{\Gamma} = \ell^\infty(\mathbbl{\Gamma})\ovot L^\infty(\check{\mathbbl{\Gamma}}).$$
   Since Morita equivalent von Neumann algebras have isomorphic centers \cite[Proposition 8.1]{Rie74} and since the center of $\ell^\infty(\mathbbl{\Gamma})$ can be identified with $\ell^\infty(\Irr(\check{\mathbbl{\Gamma}}))$,\footnote{We use the standard notation $\operatorname{Irr}(\check{\mathbbl{\Gamma}})$ for a  maximal set of pairwise non-isomorphic irreducible representations of the compact quantum group $\check{\mathbbl{\Gamma}}$.} we see that this is impossible.
   Define then 
    $$\G := \check{\mathbbl{\Gamma}}, \quad (M, \alpha):= (\ell^\infty(\mathbbl{\Gamma})\rtimes_\Delta \mathbbl{\Gamma}, \id \otimes \check{\Delta}_r), \quad (N, \beta):= (\ell^\infty(\mathbbl{\Gamma})\rtimes_\tau \mathbbl{\Gamma}, \id \otimes \check{\Delta}_r).$$
    Then by Takesaki-Takai duality, we have
    $M\rtimes_\alpha \G \sim  \ell^\infty(\mathbbl{\Gamma})\sim N\rtimes_\beta \G$. On the other hand, by Proposition \ref{crossed product}, $(M, \alpha)\not\sim_\G (N, \beta)$.
\end{Exa}

The following is the main technical result of this section. 
\begin{Prop}\label{correspondence}
    Let $\G$ be a compact quantum group and $F,G\in \operatorname{Fun}_{\Rep(\G)}(\Rep^\G(N), \Rep^\G(M))$. Then the map
    \begin{equation}\label{map}
        \operatorname{Nat}_{\Rep(\G)}(F,G)\to {}_M\mathscr{L}^\G_N(F(L^2(N)), G(L^2(N))): s \mapsto s_{L^2(N)}
    \end{equation}
    is bijective. Moreover, $s$ is a unitary natural isomorphism between $F$ and $G$ if and only if $s_{L^2(N)}$ is a unitary.
\end{Prop}
\begin{proof} Let us write
\begin{align*}
    &\mathcal{G}_F:= F(L^2(N)) = (\mathcal{G}_F, \pi_F, \rho_F, U_F) \in \Corr^\G(M,N),\\
    &\mathcal{G}_G:= G(L^2(N)) = (\mathcal{G}_G, \pi_G, \rho_G, U_G) \in \Corr^\G(M,N).
\end{align*}
    If $s,t \in \operatorname{Nat}_{\Rep(\G)}(F,G)$ with $s_{L^2(N)}= t_{L^2(N)}$, then the commutative diagram \eqref{tensor} ensures that  
    \begin{align*}
        s_{L^2(N)\otimes L^2(\G)}&= t_{L^2(N)\otimes L^2(\G)}.
    \end{align*}
    Therefore, the natural transformations $s$ and $t$ agree on a generator for $\Rep^\G(N)$, so that $s=t$. Thus, the map \eqref{map} is injective. 

   We now show that the map \eqref{map} is surjective. Given $s_1\in {}_M\mathscr{L}^\G_N(\mathcal{G}_F, \mathcal{G}_G)$, we must prove that there exists $s\in \operatorname{Nat}_{\Rep(\G)}(F,G)$ such that $s_{L^2(N)}= s_1$. We will break down the proof of the surjectivity in several smaller steps.

(STEP I) Let us start by defining 
$$s_2: = T_{L^2(N), L^2(\G)}^* \circ (s_1\otimes 1)\circ T_{L^2(N), L^2(\G)} \in {}_M\mathscr{L}^\G(F(L^2(N)\otimes L^2(\G)), G(L^2(N)\otimes L^2(\G))).$$
We claim that 
\begin{equation}\label{eqrel}
    s_2 F(U_\beta^*(\rho_N(n)\otimes 1)U_\beta) = G(U_\beta^*(\rho_N(n)\otimes 1)U_\beta) s_2, \quad n \in N.
\end{equation}
This follows from the following calculation:
\begin{align*}
    &s_2 F(U_\beta^*(\rho_N(n)\otimes 1)U_\beta)\\
    &= T_{L^2(N), L^2(\G)}^*(s_1\otimes 1) U_F^*[U_{F}T_{L^2(N), L^2(\G)} F(U_\beta^*(\rho_N(n)\otimes 1)U_\beta)T_{L^2(N), L^2(\G)}^*U_{F}^*]U_{F}T_{L^2(N), L^2(\G)}\\
    &= T_{L^2(N), L^2(\G)}^* (s_1\otimes 1)U_{F}^*(\rho_{F}(n)\otimes 1)U_{F} T_{L^2(N), L^2(\G)}\\
    &= T_{L^2(N), L^2(\G)}^* U_{G}^* (s_1 \rho_{F}(n)\otimes 1)U_{F}T_{L^2(N), L^2(\G)}\\
    &= T_{L^2(N), L^2(\G)}^* U_{G}^* (\rho_{G}(n)\otimes 1)(s_1\otimes 1)U_{F} T_{L^2(N), L^2(\G)}\\
    &= T_{L^2(N), L^2(\G)}^* U_{G}^*(\rho_{G}(n)\otimes 1)U_{G}(s_1\otimes 1) T_{L^2(N), L^2(\G)}\\
    &= T_{L^2(N), L^2(\G)}^* U_{G}^*(\rho_{G}(n)\otimes 1)U_{G} T_{L^2(N), L^2(\G)} s_2 \\
    &= G(U_\beta^*(\rho_N(n)\otimes 1)U_\beta) s_2.
\end{align*}

(STEP II) Next, define 
\begin{align*}
    s_3:&= T_{L^2(N)\otimes L^2(\G), L^2(\G)}^*(s_2\otimes 1)T_{L^2(N)\otimes L^2(\G), L^2(\G)}= T_{L^2(N), L^2(\G)\otimes L^2(\G)}^*(s_1\otimes 1 \otimes 1)T_{L^2(N), L^2(\G)\otimes L^2(\G)}
\end{align*}
which belongs to ${}_M\mathscr{L}^\G(F(L^2(N)\otimes L^2(\G)\otimes L^2(\G)), G(L^2(N)\otimes L^2(\G)\otimes L^2(\G))).$

We now claim that
\begin{equation}\label{commutation}
    s_3 F(x)= G(x) s_3, \quad x \in {}_N\mathscr{L}^\G(L^2(N)\otimes L^2(\G)\otimes L^2(\G)).
\end{equation}
From the proof of Proposition \ref{compact Morita} (4) (b), we see that
$${}_N\mathscr{L}^\G(L^2(N)\otimes L^2(\G)\otimes L^2(\G))= V_{23}^* U_{\beta,13}^* [\rho_{N\rtimes_\beta \G}(N\rtimes_\beta \G)\ovot B(L^2(\G))]_{132} U_{\beta,13}V_{23}.$$
Therefore, it suffices to show the commutation relation \eqref{commutation} for the generators
$$\begin{cases}
    V_{23}^*U_{\beta,13}^*(\rho_N(n)\otimes 1\otimes 1)U_{\beta,13}V_{23} & n \in N\\
    1\otimes V^*(1\otimes \hat{y})V & \hat{y}\in L^\infty(\hat{\mathbb{G}})\\
    1\otimes V^*(x\otimes 1)V & x \in B(L^2(\G))
\end{cases}$$
of ${}_N\mathscr{L}^\G(L^2(N)\otimes L^2(\G)\otimes L^2(\G))$. If $x\in B(L^2(\G))$ and $\hat{y}\in L^\infty(\hat{\G})$, a straightforward calculation using the fact that $V$ is a multiplicative unitary shows that 
$$V^*(x\otimes \hat{y})V\in \mathscr{L}^\G(L^2(\G)\otimes L^2(\G)).$$
Therefore, we find
\begin{align*}
    s_3 F(1\otimes V^*(x\otimes \hat{y})V) &= T_{L^2(N), L^2(\G)\otimes L^2(\G)}^*(s_1\otimes 1 \otimes 1)T_{L^2(N), L^2(\G)\otimes L^2(\G)} F(1\otimes V^*(x\otimes \hat{y})V)\\
    &=T_{L^2(N), L^2(\G)\otimes L^2(\G)}^*(s_1\otimes 1\otimes 1) (1\otimes V^*(x\otimes \hat{y})V) T_{L^2(N), L^2(\G)\otimes L^2(\G)}\\
    &= T_{L^2(N), L^2(\G)\otimes L^2(\G)}^* (1\otimes V^*(x\otimes \hat{y})V) (s_1\otimes 1\otimes 1) T_{L^2(N), L^2(\G)\otimes L^2(\G)}\\
    &= G(1\otimes V^*(x\otimes \hat{y})V) T_{L^2(N), L^2(\G)\otimes L^2(\G)}^* (s_1\otimes 1 \otimes 1)T_{L^2(N), L^2(\G)\otimes L^2(\G)}\\
    &= G(1\otimes V^*(x\otimes \hat{y})V) s_3.
\end{align*}
It remains to prove the commutation relation 
\begin{equation}
    s_3 F(V_{23}^* U_{\beta,13}^*(\rho_N(n)\otimes 1\otimes 1)U_{\beta,13}V_{23}) = G(V_{23}^* U_{\beta,13}^*(\rho_N(n)\otimes 1\otimes 1)U_{\beta,13}V_{23}) s_3, \quad n \in N.
\end{equation}
As before, it is convenient to simplify some notation. Therefore, we write 
\begin{align*}
    T_V&:= T_{L^2(N), L^2(\G)}= T_{L^2(N), (L^2(\G),V)}\\
    T_{V_{13}V_{23}}&:= T_{L^2(N), L^2(\G)\otimes L^2(\G)}= T_{L^2(N), (L^2(\G)\otimes L^2(\G), V_{13}V_{23})}\\
    T_{V_{13}}&:= T_{L^2(N), (L^2(\G)\otimes L^2(\G), V_{13})}\\
    T_{\mathbb{I}}^{\otimes}&:= T_{L^2(N)\otimes L^2(\G), (L^2(\G), \mathbb{I})}.
\end{align*}

A straightforward calculation shows that, with $\Sigma: L^2(\G)\otimes L^2(\G)\to L^2(\G)\otimes L^2(\G): \xi\otimes \eta \mapsto \eta \otimes \xi,$ $$V^* \Sigma \in \mathscr{L}^\G((L^2(\G)\otimes L^2(\G), V_{13}), (L^2(\G)\otimes L^2(\G), V_{13}V_{23})).$$
Therefore, making multiple times use of naturality, and keeping in mind \eqref{eqrel}, we find:
\begin{align*}
    &s_3 F(V_{23}^* U_{\beta,13}^* (\rho_N(n)\otimes 1 \otimes 1)U_{\beta,13}V_{23})\\
    &= T_{V_{13}V_{23}}^*(s_1\otimes 1 \otimes 1)T_{V_{13}V_{23}} F(1\otimes V^* \Sigma) F(U_{\beta,12}^*(\rho_N(n)\otimes 1\otimes 1)U_{\beta,12}) F(1\otimes \Sigma V) \\
    &= T_{V_{13}V_{23}}^* (s_1\otimes 1 \otimes 1)(1\otimes V^*\Sigma)T_{V_{13}} F((U_\beta^*(\rho_N(n)\otimes 1)U_\beta) \otimes 1) F(1\otimes \Sigma V)\\
    &= T_{V_{13}V_{23}}^* (1\otimes V^*\Sigma)(s_1\otimes 1 \otimes 1)(T_V\otimes 1)T_{\mathbb{I}}^{\otimes} F((U_\beta^*(\rho_N(n)\otimes 1)U_\beta) \otimes 1) F(1\otimes \Sigma V)\\
    &= T_{V_{13}V_{23}}^* (1\otimes V^*\Sigma)(s_1\otimes 1 \otimes 1)(T_V\otimes 1)(F(U_\beta^*(\rho_N(n)\otimes 1)U_\beta) \otimes 1)T_{\mathbb{I}}^{\otimes} F(1\otimes \Sigma V)\\
    &= T_{V_{13}V_{23}}^* (1\otimes V^*\Sigma) (T_V\otimes 1)(s_2\otimes 1)(F(U_\beta^*(\rho_N(n)\otimes 1)U_\beta)\otimes 1)T_{\mathbb{I}}^\otimes F(1\otimes \Sigma V)\\
    &= T_{V_{13}V_{23}}^* (1\otimes V^*\Sigma)(T_V\otimes 1) (G(U_\beta^*(\rho_N(n)\otimes 1)U_\beta)\otimes 1)(s_2\otimes 1)T_{\mathbb{I}}^\otimes F(1\otimes \Sigma V)\\
    &= T_{V_{13}V_{23}}^*(1\otimes V^*\Sigma) T_{V_{13}}T_{\mathbb{I}}^{\otimes*} (G(U_\beta^*(\rho_N(n)\otimes 1)U_\beta)\otimes 1)(s_2\otimes 1)T_{\mathbb{I}}^\otimes F(1\otimes \Sigma V)\\
    &= T_{V_{13}V_{23}}^*(1\otimes V^*\Sigma)T_{V_{13}} G(U_{\beta,12}^*(\rho_N(n)\otimes 1 \otimes 1)U_{\beta,12})T_{\mathbb{I}}^{\otimes *} (s_2\otimes 1)T_{\mathbb{I}}^\otimes F(1\otimes \Sigma V)\\
    &= T_{V_{13}V_{23}}^*(1\otimes V^*\Sigma)T_{V_{13}} G(U_{\beta,12}^*(\rho_N(n)\otimes 1 \otimes 1)U_{\beta,12})T_{\mathbb{I}}^{\otimes *} (s_2\otimes 1) (T_{V}^*\otimes 1)T_{V_{13}} F(1\otimes \Sigma V)\\
 &= T_{V_{13}V_{23}}^*(1\otimes V^*\Sigma)T_{V_{13}} G(U_{\beta,12}^*(\rho_N(n)\otimes 1 \otimes 1)U_{\beta,12})T_{\mathbb{I}}^{\otimes *} (T_V^*\otimes 1) (s_1\otimes 1 \otimes 1)T_{V_{13}} F(1\otimes \Sigma V)\\
 &= G(1\otimes V^*\Sigma) G(U_{\beta,12}^*(\rho_N(n)\otimes 1 \otimes 1)U_{\beta,12}) T_{V_{13}}^*(s_1\otimes 1 \otimes 1) T_{V_{13}} F(1\otimes \Sigma V)\\
 &= G(1\otimes V^*\Sigma) G(U_{\beta,12}^*(\rho_N(n)\otimes 1 \otimes 1)U_{\beta,12}) T_{V_{13}}^*(s_1\otimes 1 \otimes 1) (1\otimes \Sigma V)T_{V_{13}V_{23}} \\
 &= G(1\otimes V^*\Sigma) G(U_{\beta,12}^*(\rho_N(n)\otimes 1 \otimes 1)U_{\beta,12}) G(1\otimes \Sigma V)T_{V_{13}V_{23}}^*(s_1\otimes 1 \otimes 1)T_{V_{13}V_{23}} \\
 &= G(V_{23}^*U_{\beta,13}^*(\rho_N(n)\otimes 1 \otimes 1)U_{\beta, 13}V_{23})s_3.
\end{align*}

(STEP III) By Proposition \ref{compact Morita} (4) (b), $L^2(N)\otimes L^2(\G)\otimes L^2(\G)$ is a generator for $\Rep^\G(N)$. Therefore, Proposition \ref{compact Morita} (5) and \eqref{commutation} imply that there exists a unique natural transformation $s\in \operatorname{Nat}(F,G)$ such that $s_3=s_{L^2(N)\otimes L^2(\G)\otimes L^2(\G)}$. We now argue that $s_1 = s_{L^2(N)}$.

Given an equivariant representation $(\mathcal{H},U)$, we define the isometric $\G$-intertwiner
$$\iota_{\mathcal{H}}: (\mathcal{H}, U)\to (\mathcal{H}\otimes L^2(\G), U_{13}V_{23}): \xi \mapsto \xi \otimes \xi_\G.$$
 Using the fact that the isomorphisms $T_{\mathcal{H}, \mathcal{K}}$ are natural and using the commutativity of the diagram \eqref{commutativity2}, we easily deduce that the diagram
\begin{equation}\label{dia}
    \begin{tikzcd}
F(\mathcal{H}) \arrow[rd, "\iota_{F(\mathcal{H})}"'] \arrow[rr, "F(\iota_H)"] &                               & F(\mathcal{H}\otimes L^2(\G)) \arrow[ld, "{T_{\mathcal{H}, L^2(\G)}}"] \\
                                                                              & F(\mathcal{H})\otimes L^2(\G) &                                                                       
\end{tikzcd}
\end{equation}
commutes (and similarly for the functor $G$). We then compute 
\begin{align*}
    s_{L^2(N)\otimes L^2(\G)}&= G(\iota_{L^2(N)\otimes L^2(\G)})^* \circ s_{L^2(N)\otimes L^2(\G)\otimes L^2(\G)}\circ F(\iota_{L^2(N)\otimes L^2(\G)})\\
    &= G(\iota_{L^2(N)\otimes L^2(\G)})^* \circ T_{L^2(N)\otimes L^2(\G), L^2(\G)}^* \circ (s_2\otimes \id) \circ T_{L^2(N)\otimes L^2(\G), L^2(\G)}\circ  F(\iota_{L^2(N)\otimes L^2(\G)})\\
    &= \iota_{G(L^2(N)\otimes L^2(\G))}^*\circ (s_2\otimes \id) \circ  \iota_{F(L^2(N)\otimes L^2(\G))} = s_2
\end{align*}
where the first equality follows from naturality of $s$ and the fact that $\iota_{L^2(N)\otimes L^2(\G)}$ is isometric, the second equality from the definition of $s_3$ and the third equality from the commutativity of the diagram \eqref{dia}.

Completely similarly, we show that $s_{L^2(N)}= s_1$. 

(STEP IV) We now argue that the diagram \eqref{tensor} commutes, so that $s\in \operatorname{Nat}_{\Rep(\G)}(F,G).$
By construction, the diagram \eqref{tensor} commutes if $\mathcal{H}= L^2(N)\otimes L^2(\G)\in \Rep^\G(N)$ and $\mathcal{K}= L^2(\G)\in \Rep(\G)$. Therefore, \eqref{tensor} also commutes for $\mathcal{H}$ an arbitrary direct sum of $L^2(N)\otimes L^2(\G)$ and $\mathcal{K}$ an arbitrary direct sum of $L^2(\G)$ (here, we implicitly use that $F,G$ preserve arbitrary direct sums, cfr.\ \cite[Proposition 4.9]{Rie74}). If now $\mathcal{H}\in \Rep^\G(N)$ and $\mathcal{K}\in \Rep(\G)$ are arbitrary, we can find isometries 
$$v\in {}_N\mathscr{L}^\G(\mathcal{H}, \mathcal{H'}), \quad w\in \mathscr{L}^\G(\mathcal{K},\mathcal{K'})$$
such that \eqref{tensor} commutes for $\mathcal{H'}\in \Rep^\G(N)$ and $\mathcal{K'}\in \Rep(\G)$. Consider then the following diagram:
$$
\begin{tikzcd}
F(\mathcal{H}\otimes \mathcal{K}) \arrow[rd, "F(v\otimes w)"] \arrow[ddd, "s_{\mathcal{H}\otimes \mathcal{K}}"', bend right] \arrow[rrrr, "{T_{\mathcal{H}, \mathcal{K}}}"] &                                                                                                                                      &  &                                                                              & F(\mathcal{H})\otimes \mathcal{K} \arrow[ld, "F(v)\otimes w"'] \arrow[ddd, "s_{\mathcal{H}}\otimes \id", bend left] \\
                                                                                                                                                                            & F(\mathcal{H}'\otimes \mathcal{K'}) \arrow[rr, "{T_{\mathcal{H'}, \mathcal{K'}}}"] \arrow[d, "s_{\mathcal{H'}\otimes\mathcal{K'}}"'] &  & F(\mathcal{H}')\otimes \mathcal{K}' \arrow[d, "s_{\mathcal{H'}}\otimes \id"] &                                                                                                                     \\
                                                                                                                                                                            & G(\mathcal{H}'\otimes \mathcal{K}') \arrow[rr, "{T_{\mathcal{H'}, \mathcal{K'}}}"]                                                   &  & G(\mathcal{H}')\otimes \mathcal{K'}                                          &                                                                                                                     \\
G(\mathcal{H}\otimes \mathcal{K}) \arrow[ru, "G(v\otimes w)"'] \arrow[rrrr, "{T_{\mathcal{H}, \mathcal{K}}}"']                                                              &                                                                                                                                      &  &                                                                              & G(\mathcal{H})\otimes \mathcal{K} \arrow[lu, "G(v)\otimes w"]                                                      
\end{tikzcd}$$
Our goal is to show that the outer diagram commutes. We know that the inner rectangle commutes. By naturality, the four outer subdiagrams commute as well. A diagram chase then yields
$$(G(v)\otimes w)\circ T_{\mathcal{H}, \mathcal{K}}\circ s_{\mathcal{H}\otimes \mathcal{K}}= (G(v)\otimes w) \circ (s_{\mathcal{H}}\otimes \id)\circ T_{\mathcal{H}, \mathcal{K}}$$
which implies that the outer diagram commutes, since $v$ and $w$ are isometries. We conclude that $s\in \operatorname{Nat}_{\Rep(\G)}(F,G)$. This proves the surjectivity of the map \eqref{map}.

    Finally, assume that $s_1\in {}_M\mathscr{L}_N^\G(F(L^2(N)), G(L^2(N)))$, and let $s\in \operatorname{Nat}_{\Rep(\G)}(F,G)$ be such that $s_1= s_{L^2(N)}$. Consider $s^*\in \operatorname{Nat}_{\Rep(\G)}(G,F)$ defined by
    $$(s^*)_{\mathcal{H}}:= s_{\mathcal{H}}^*\in {}_M\mathscr{L}^\G(G(\mathcal{H}), F(\mathcal{H})), \quad \mathcal{H}\in \Rep^\G(N).$$
    Then we have $(s^* \circ s)_{L^2(N)}= \id_{F(L^2(N))}$ and $(s\circ s^*)_{L^2(N)} = \id_{G(L^2(N))}$ so that $s^* \circ s = \id_{F}$ and $s\circ s^*= \id_{G}$ by the injectivity of \eqref{map}. This shows that $s$ is a unitary isomorphism.
\end{proof}

The following result is then the announced equivariant version of the Eilenberg-Watts theorem:

\begin{Theorem}\label{main}
    Let $\G$ be a compact quantum group. The functors
    \begin{align*}
        &P: \operatorname{Corr}^\G(M,N)\to \operatorname{Fun}_{\Rep(\G)}(\Rep^\G(N), \Rep^\G(M)): \mathcal{G}\mapsto F_{\mathcal{G}}, \\
    &Q: \operatorname{Fun}_{\Rep(\G)}(\Rep^\G(N), \Rep^\G(M))\to  \operatorname{Corr}^\G(M,N): F \mapsto F(L^2(N))
    \end{align*}
    are quasi-inverse to each other. 
    
    In particular, if $F: \Rep^\G(N)\to \Rep^\G(M)$ a $\Rep(\G)$-module $*$-functor and $\mathcal{G}:= F(L^2(N))\in \Corr^\G(M,N)$, there is a unique unitary natural isomorphism $s_F\in \operatorname{Nat}_{\Rep(\G)}(F, F_{\mathcal{G}})$ extending the natural isomorphism $F(L^2(N)) = \mathcal{G}\cong \mathcal{G}\boxtimes_N L^2(N) = F_{\mathcal{G}}(L^2(N))$.
\end{Theorem}
\begin{proof} We only need to prove that $PQ\cong \id$. Consider a $\Rep(\G)$-module $*$-functor $F: \Rep^\G(N)\to \Rep^\G(M)$. Consider the isomorphism
$$s_{F,1}: [PQ(F)](L^2(N))= F_{Q(F)}(L^2(N)) = Q(F)\boxtimes_N L^2(N)\cong Q(F) = F(L^2(N)): y \otimes_N  \eta \mapsto y\eta$$
of $\G$-$M$-$N$-correspondences. By Proposition \ref{correspondence}, we deduce that there is a unique unitary isomorphism $s_F \in \operatorname{Nat}_{\Rep(\G)}(PQ(F), F)$ such that $(s_F)_{L^2(N)}= s_{F,1}$. 

It remains to show that these unitary isomorphisms are natural in $F\in \operatorname{Fun}_{\Rep(\G)}(\Rep^\G(N), \Rep^\G(M))$, i.e. we must show that if $s\in \operatorname{Nat}_{\Rep(\G)}(F,G)$, then the diagram
\begin{equation}
    \label{commutative diagram}
\begin{tikzcd}
PQ(F) \arrow[rr, "s_F"] \arrow[d, "PQ(s)"'] &  & F \arrow[d, "s"] \\
PQ(G) \arrow[rr, "s_G"']                    &  & G               
\end{tikzcd}
\end{equation}

commutes. Note however that for $y\in \mathscr{L}_N(L^2(N), Q(F))$ and $\eta \in L^2(N)$, we have
\begin{align*}
    s_{L^2(N)} s_{F,1} (y\otimes_N \eta) &= s_{L^2(N)} (y\eta) = s_{G,1} (s_{L^2(N)}y \otimes_N \eta)= s_{G,1} (PQ)(s)_{L^2(N)} (y\otimes_N \eta)
\end{align*}
so by injectivity of \eqref{map}, we see that the diagram \eqref{commutative diagram} commutes. 
\end{proof}

\textbf{Acknowledgments:} The research of the author was supported by Fonds voor Wetenschappelijk Onderzoek (Flanders), via an FWO Aspirant fellowship, grant 1162524N. The author would like to thank K. De Commer for useful input and discussions
throughout this entire project. We thank J. Krajczok, L. Rollier and J. Vercruysse for a useful discussion.

\end{document}